\documentclass{amsart}

\usepackage{amssymb}

\makeatletter
\renewcommand{\subsection}{
\@startsection%
{subsection}
{2}
{0em}
{2ex plus 0.1ex minus -0.05ex}
{-1em plus 0.2em}
{
\bfseries}}
\makeatother

\makeatletter                                            
\renewcommand{\@begintheorem}[2]{                        
\rm \trivlist \item [\hskip \labelsep {\bf #2\ \ #1.}]   
                                }                        
\makeatother                                             

\newcommand{\ts}{\vspace{\baselineskip}\noindent{\bf Proof.}$\;\;$}

\newcommand{\ZZ}{{\mathbb Z}}
\newcommand{\QQ}{{\mathbb Q}}
\newcommand{\RR}{{\mathbb R}}

\newcommand{\CC}{{\mathbb C}}

\newcommand{\PP}{{\bf P}}

\newcommand{\cE}{{\mathcal E}}

\newcommand{\cI}{{\mathcal I}}

\newcommand{\cL}{{\mathcal L}}
\newcommand{\cM}{{\mathcal M}}

\newcommand{\cO}{{\mathcal O}}

\newcommand{\cT}{{\mathcal T}}

\newcommand{\om}{\omega}

\newcommand{\Hom}{{\operatorname{Hom}}}
\newcommand{\Nm}{\mbox{Nm}}

\begin{document}

\title{Fourfolds of Weil type and the spinor map}
\author{Bert van Geemen}
\address{Dipartimento di Matematica, Universit\`a di Milano,
Via Saldini 50, I-20133 Milano, Italia}
\email{lambertus.vangeemen@unimi.it}

\begin{abstract}
Recent papers by Markman and O'Grady give,
besides their main results on the Hodge conjecture and 
on hyperk\"ahler varieties,
surprising and explicit
descriptions of families of abelian fourfolds of Weil type with trivial 
discriminant.
They also provide a new perspective on the well-known 
fact that these abelian varieties are Kuga Satake varieties for certain weight 
two Hodge structures of rank six. 

In this paper we give a pedestrian introduction
to these results. The spinor map, which is defined using a half-spin
representation of $SO(8)$, is used intensively. 
For simplicity, we use basic representation theory and we avoid the use 
of triality.
\end{abstract}

\maketitle

\section*{Introduction}
The recent papers \cite{Markman}, \cite{O'G} by Markman and O'Grady
provide new descriptions of families of abelian fourfolds of Weil type. 
Markman uses these to prove that certain Hodge classes on these fourfolds are algebraic.
Both show that these abelian varieties are isogeneous to the intermediate
Jacobians of algebraic hyperk\"ahler varieties of Kummer type.
O'Grady further relates this to the 
Kuga Satake construction for the (primitive) second cohomology group 
of algebraic Kummer type varieties. See also \cite{Voisin} for further developments.

An abelian fourfold of Weil type has an imaginary quadratic field 
$K=\QQ(\sqrt{-d})$ in its endomorphism algebra.
These fourfolds define two subspaces of
the complexification of their first homology group $H_1$, 
a free $\ZZ$-module of rank $8$.
They are the $+i$-eigenspace of the complex structure on $H_1\otimes\RR$ defined by $A$
and one of the two eigenspaces of the $K$-action. 
Markman obtains the polarization on the abelian fourfold,
an alternating form on $H_1$, from a symmetric(!) form on $H_1$
and the $K$-action. The two subspaces turn out to be maximally isotropic
subspaces for this symmetric form.

In this paper we will mainly follow Markman's approach. 
He considers a free, rank $8$, $\ZZ$-module $V$ 
equipped with a bilinear form. This $V$ will become the first 
cohomology group of the fourfolds of Weil type. The maximally
isotropic subspaces of the complexification $V_\CC$ of $V$ are well-known
to be parametrized by two copies of a Legendrian Grassmannian, a complex
manifold of dimension six. The spinor map is a natural embedding of this Grassmannian in $\PP^7$, 
the image is a quadric hypersurface $Q^+$. 
This map already made several appearances in algebraic geometry,
for example in the study of vector bundles over hyperelliptic curves in
\cite{vG_VB}, of K3 surfaces in \cite{MukaiCSP1}, of secant varieties
in \cite{ManivelSpinor} and of integrable systems \cite{BaloghHH}.

The spinor map is best constructed using the representation theory
of $Spin(V)$, a double cover of the orthogonal group $SO(V)$
defined by the bilinear form on $V$. 
The spin group has a half-spin representation whose projectivization is $\PP^7$.
The spinor map is equivariant for the action of $Spin(V)$.
A natural integral structure on the half-spin representation
allows one to identify it with the complexification
of a free $\ZZ$ module $S^+$ of rank $8$. 
There is a non-degenerate bilinear form on $S^+$ which defines the quadric $Q^+$.

An analytic open subset $\Omega\subset Q^+$ parametrizes complex structures on $V_\RR$
that preserve the bilinear form on $V$.
Fixing a general element $s\in S^+$ and considering only the complex structures 
on $V_\RR$ corresponding to $\ell\in \Omega\cap s^\perp$ produces
a five dimensional family of complex tori $\cT_\ell$, not algebraic in general, 
that have a Hodge class, called the Cayley class, 
$$
c_s\,\in\, H^{2,2}(\cT_\ell,\ZZ)\,=\,H^4(\cT_\ell,\ZZ) \cap H^{2,2}(\cT_\ell).
$$
The idea of using these tori and the associated action of $Spin(7)=Spin(s^\perp)$
to study the Hodge conjecture for fourfolds of Weil type is due to V.\ Mu\~noz \cite{Munoz}.
In \S \ref{cayleyclasses} we observe that the existence of the Cayley classes 
can be deduced from a relation between the spinor and the Pl\"ucker map. 
Using representation theory we then compute the class $c_s$ for certain $s$ that 
are relevant for Markman's results in Proposition \ref{explicitCs}.

For any $h\in S^+$ such that the sublattice $\langle h,s\rangle$ of $S^+$ spanned by $h$ and $s$ 
has rank two and is positive definite, 
the tori parametrized by the four dimensional domain 
$\Omega\cap \langle h,s\rangle ^\perp$ turn out to be abelian fourfolds of Weil type. 
The imaginary quadratic field $K$ depends on $h,s$, but fixing $s$ and
choosing $h$ suitably, any such field occurs. 
The polarization is determined by $K$ and the bilinear form on $V$.
A further discrete invariant, the discriminant of a polarized abelian variety of Weil type, 
is always trivial for the fourfolds constructed in this way. 
See Theorem \ref{Weil} for these results of Markman and O'Grady.

The Hodge conjecture for an abelian fourfold $A$ of Weil type is non-trivial. 
There is a natural 2-dimensional
subspace $W_K\subset H^{2,2}(A,\QQ)$ of Hodge classes.
It is not known in general if this subspace is spanned by classes of algebraic
cycles. If $c_s$ is algebraic, then all classes in $W_K$ are also algebraic. 
Markman makes important progress in the study of the Hodge conjecture by
showing that $c_s$ is algebraic for all abelian fourfolds appearing in his construction, 
which are all fourfolds of Weil type with trivial discriminant.
For this he uses deformation theory of sheaves on hyperk\"ahler manifolds, 
see \S \ref{moduli} for a brief outline.

Triality, an automorphism of order three of $Spin(V)$, allows one to
relate the standard representation of $Spin(V)$ (via $(SO(V)$ on $V$)
and the two half-spin representations, one of which is $S^+$. 
While it is prominent in \cite{Markman}, we use instead an `ad hoc' Lemma \ref{acth}.
It is of importance for instance in the results on the Cayley class and for the Kuga Satake varieties.

We limit ourselves to  a basic exposition of the constructions of Markman and O'Grady 
of the abelian fourfolds of Weil type with trivial discriminant and of the Cayley classes
of Mu\~noz and Markman. 
Some details of the representation theory involved 
in the construction of the spinor map and the Cayley classes
can be found in the Appendix, \S \ref{App}.
The relation with the Kuga Satake construction is indicated in \S \ref{KSvar} - \ref{H2ell}.

\subsection*{Acknowledgements} Discussions with E.\ Markman and K.G.\ O'Grady 
were very helpful.

\quad

\section{Tori with an orthogonal structure}\label{spinorweil}

\subsection{The lattice $V$}\label{defV}
The complex tori we consider are all quotients of a fixed real vector space, 
with a varying complex structure, by a fixed lattice. 
Whereas one might expect an alternating form, a polarization, on the first cohomology group
to be important, Markman instead fixes a symmetric, non-degenerate, bilinear form 
$(\bullet,\bullet)_V$ on a rank eight free $\ZZ$-module $V$ of signature $(4+,4-)$.
He fixes a rank four free $\ZZ$-module $W$, defines 
$W^*:=\Hom_\ZZ(W,\ZZ)$ and 
$$
V\,:=\,W\,\oplus\, W^*,\qquad \big((w_1,w_1^*),(w_2,w_2^*)\big)_V\,:=\,w_1^*(w_2)\,+\,w_2^*(w_1)~.
$$
If $e_1,\ldots,e_4$ is a $\ZZ$-basis of $W$ and $e_{i+4}:=e_i^*$, where $e_1^*,\ldots,e_4^*$ 
is the dual basis of $W^*$
so that $e_i^*(e_j)=\delta_{ij}$ (Kronecker's delta), then 
$$
(v_1,v_2)_V\,:=\,\sum_{i=1}^4x_iy_{i+4}+x_{i+4}y_i,\qquad
\big(v_1:=\sum_{i=1}^8 x_ie_i,\quad v_2:=\sum_{i=1}^8 y_ie_i\;\in V\big)~,
$$
hence $(V,(\bullet,\bullet)_V)\cong U^{\oplus 4}$, the direct sum of four copies of the 
hyperbolic plane $U=(\ZZ^2,\small{\begin{pmatrix}0&1\\1&0\end{pmatrix}})$.

In \cite{Markman} one finds $W:=H^1(X,\ZZ)$ for an abelian surface $X$, 
but for the basic properties of the complex tori this is not needed.

\subsection{Complex structures on $V_\RR$} Let $V_\RR:=V\otimes_\ZZ\RR$, 
it is an eight dimensional vector space over the real numbers. 
A complex structure on $V_\RR$ is a linear map $J:V_\RR\rightarrow V_\RR$ with $J^2=-I$.
Such a map has two (complex) eigenspaces $Z_+,Z_-\subset V_\CC:=V\otimes_\ZZ\CC$ corresponding
to the eigenvalues $i,-i\in\CC$ of $J$. These eigenspaces are complex conjugate, 
$\overline{Z_+}=Z_-$, where the complex conjugation on $V_\CC$
is defined as $\overline{v\otimes z}=v\otimes \bar{z}$ for $v\in V$ and $z\in\CC$.
$$
V_\CC\,=\,Z_+\oplus Z_-=Z_+\oplus\overline{Z_+},\qquad 
J\,=\,(i,-i)\,\in\,End(Z_+)\oplus End(Z_-)~.
$$

Conversely, given two complex conjugate subspaces $Z_\pm\subset V_\CC$ such that 
$V_\CC=Z_+\oplus Z_-$
one can define a linear map $\tilde{J}:V_\CC\rightarrow V_\CC$ by $\tilde{J}(v_++v_-)=
iv_+-iv_-$. Then there is a linear map $J:V_\RR\rightarrow V_\RR$ whose $\CC$-linear extension
to $V_\CC$ is $\tilde{J}$. In fact, the inclusion $V_\RR\hookrightarrow V_\CC$ identifies
$V_\RR$ with the $(v_+,v_-)\in Z_+\oplus Z_-$ with $\overline{v_+}=v_-$.  
Writing $v\in V_\RR$ as $v=v_+ +v_-$, with $v_-=\overline{v_+}$, one has
$\tilde{J}v=iv_++\overline{iv_+}\in V_\RR$, so $J$ is just the restriction of $\tilde{J}$ to $V_\RR$.

\subsection{Orthogonal complex structures and isotropic subspaces}
The $\RR$-bilinear extension of $(\bullet,\bullet)_V$ defines
a bilinear form on $V_\RR$, denoted by the same symbol. 
We consider now the complex structures $J$
that preserve this bilinear form, so
$(Jv_1,Jv_2)_V=(v_1,v_2)_V$ for all $v_1,v_2\in V_\RR$. 
Equivalently, $J\in SO(V_\RR,(\bullet,\bullet)_V)$ and we will call $J$ an orthogonal complex structure. 
Notice that for such a complex structure $J$ and for eigenvectors $v_{1+},v_{2,+}\in Z_+$ we have,
for the $\CC$-bilinear extension of the bilinear form,
$$
(v_{1+},v_{2+})_V\,=\,
(Jv_{1+},Jv_{2+})_V\,=\,(iv_{1+},iv_{2+})_V\,=\,i^2(v_{1+},v_{2+})_V\,=\,-(v_{1+},v_{2+})_V~.
$$
Hence the restriction of $(\bullet,\bullet)_V$ to $Z_+$ is identically zero. Thus $Z_+$
is an isotropic subspace of $V_\CC$ (and since $\dim Z_+=4=(1/2)\dim V_\CC$ it is a maximally
isotropic, or Legendrian, subspace of $V_{\CC}$). Similarly $Z_-$ is a maximally isotropic subspace
of $V_\CC$ (and since the bilinear form is non-degenerate it induces a duality $Z_+\cong Z_-^*$). 

One easily verifies that, conversely,
an isotropic subspace $Z_+\subset V_\CC$ such that $V_\CC=V_+\oplus\overline{V_+}$
defines a complex structure $J$ on $V_\RR$ that preserves $(\bullet,\bullet)_V$.
We summarize this in the following lemma.

\subsection{Lemma}
There is a natural bijection between the following two sets:
\begin{itemize}
 \item the orthogonal complex structures $J\in SO(V_\RR,(\bullet,\bullet)_V)$
on $V_\RR$,
\item the maximally isotropic subspaces 
$Z$ of $V_\CC$ such that $V_\CC=Z\oplus\overline{Z}$.
\end{itemize}

\

\section{The Legendrian Grassmannian and the spinor map}

\subsection{}  In this section we recall that a connected component $IG(4,V_\CC)^+$
of the Grassmannian 
of maximally isotropic subspaces of $V_\CC$
is isomorphic to a smooth six dimensional quadric $Q^+\subset \PP S ^+_\CC\cong \PP^7$,
where $(S^+,(\bullet,\bullet)_{S^+})$ is a certain lattice of rank eight. 
This isomorphism is induced by the spinor map, which is equivariant for the action
of the double cover $Spin(V)$ of $SO(V)$ on $V$ and $S^+$ respectively.
We refer to the Appendix \S \ref{App} for more details.

\subsection{The Grassmannian $IG(4,V_\CC)^+$}
The (complex) four dimensional subspaces of $V_\CC$ are parametrized by the Grassmannian 
$G(4,V_\CC)$,
which has dimension $4\cdot(8-4)=16$.
The maximally isotropic subspaces for $(\bullet,\bullet)_V$ 
(which are also those for the associated quadratic form) are parametrized by two
(isomorphic, disjoint, connected) complex submanifolds of dimension six of $G(4,V_\CC)$, 
denoted by $IG(4,V_\CC)^+$ and $IG(4,V_\CC)^-$. 
(See \cite[Chapter 6]{GH} for linear subspaces of quadrics.)
This generalizes the two rulings 
(families of lines) on a smooth quadric $Q\cong\PP^1\times\PP^1$ in $\PP^3$.
We denote by $IG(4,V_\CC)^+$ 
the connected component which contains the maximally isotropic subspace $W_\CC^*$.
A complex maximally isotropic subspace $Z$ defines a point $[Z]\in IG(4,V_\CC)^+$ 
if and only if the dimension of $Z\cap W_\CC^*$ is even.
In particular, also $[W_\CC]\in IG(4,V_\CC)^+$.

We recall a local parametrization of $IG(4,V_\CC)^+$ by 
alternating $4\times 4$ complex matrices. 
A basis of $W^*$ is given by the last four basis vectors of $V$ in \S \ref{defV}. Thus $W_\CC^*$
is spanned by the columns of the $8\times 4$ matrix $({}^0_I)$. 
Slightly deforming $W_\CC$,
we obtain another subspace spanned by the columns of an $8\times 4$ matrix. 
Since $det I=1\neq 0$ we may assume that the lower 
$4\times 4$ submatrix is still invertible. Then we can find a basis of the same subspace 
given by the columns of a matrix $({}^B_I)$, the corresponding subspace will be denoted by $Z_B$. 
Thus we found a Zariski open subset $G(4,V_\CC)_0$ of $G(4,V_\CC)$ of dimension $4^2=16$
parametrized by $4\times 4$ complex matrices.

In general $Z_B$ will not be isotropic, but one easily verifies that 
$$
{(\bullet,\bullet)_V}_{|Z_B\times Z_B}\,=\,0 
\quad \Longleftrightarrow\quad
\begin{pmatrix}^tB&I\end{pmatrix}\begin{pmatrix}0&I\\I&0\end{pmatrix}
\begin{pmatrix}B\\I\end{pmatrix}\,=\,0
\quad \Longleftrightarrow\quad{}^tB\,+\,B\,=\,0~.
$$ 
Hence the vector space of alternating $4\times 4$ matrices $Alt_4$ 
provides us with a parametrization of a Zariski open subset of $IG(4,V_\CC)^+$ 
of dimension $4(4-1)/2=6$ which we denote by $IG(4,V_\CC)^+_0$:
$$
Alt_4\,\stackrel{\cong}{\longrightarrow}\,IG(4,V_\CC)^+_0\;\hookrightarrow\;IG(4,V_\CC)^+,
\qquad B\,\longmapsto\, [Z_B]\,=\,\left[\begin{pmatrix}B\\I\end{pmatrix}\right]~.
$$

The isotropic subspace $Z_B$ is also the graph of the (alternating) map
$W^*\rightarrow W$, $w^*\mapsto B w^*$.
 
\

\subsection{The Pl\"ucker map}\label{plueckermap} 
The Grassmannian $G(4,V_\CC)$ has a natural embedding, 
the Pl\"ucker map $\pi$,
into a projective space $\PP^N=\PP\wedge^4V_\CC$ of dimension $N+1=\binom{8}{4}=70$:
$$
\pi:\,G(4,V_\CC)\,\longrightarrow\,\PP \wedge^4V_\CC,\qquad 
Z\,\longmapsto \,[\wedge^4Z]~.
$$
The Pl\"ucker map is equivariant for the action of $GL(V_\CC)$.

On the open subset $G(4,V_\CC)_0$ of $G(4,V_\CC)$ the Pl\"ucker map is thus given
by the determinants of the $4\times 4$ submatrices of the $8\times 4$ matrix $P:=(^B_I)$. 
Using the basis of $V$ from \S \ref{defV}, 
the coefficient of $e_{i_1}\wedge\ldots\wedge e_{i_4}$ in 
$$
[\wedge^4 Z_B]\,=\,[r_1\wedge\ldots\wedge r_4],\qquad 
\big(r_j\,=\,\sum_{k=1}^8P_{kj}e_k\,\in\,Z_B,\quad P:=(^B_I)\big)
$$
is the determinant of the $4\times 4$ submatrix of $P$ with rows $i_1,\ldots,i_4$.

\

\subsection{The spinor map}\label{spinormap} 
The Picard group of $G(4,V_\CC)$ is generated by the Pl\"ucker line bundle $\pi^*\cO_{\PP^N}(1)$.
The restriction of this line bundle to $IG(4,V_\CC)^+$ does not generate the Picard group of
$IG(4,V_\CC)^+$, but there is a line bundle $\cL$ on $IG(4,V_\CC)^+$ such that 
$$
(\pi^*\cO_{\PP^N}(1))_{|IG(4,V_\CC)^+}\,\cong\,\cL^{\otimes 2}~,
$$
and $Pic(IG(4,V_\CC)^+)\cong\ZZ$ is generated by $\cL$. One has $h^0(IG(4,V_\CC)^+,\cL)=8$ and
the natural spinor (or Cartan) map 
$$
\gamma:\,IG(4,V_\CC)^+\,\longrightarrow\,\PP S^+_\CC\,\cong\,\PP H^0(IG(4,V_\CC)^+,\cL)^*
$$
is an embedding whose image is a smooth quadric $Q^+\subset\PP S^+_\CC$. 
Here $S^+_\CC$ is the complexification of a lattice $S^+$ 
that will be defined below (in an ad hoc manner),
see also \S \ref{MukaiLattice} and \S \ref{cliffordV}.

For $Z_B$ in the open subset $IG(4,V_\CC)^+_0$, where $B=(b_{ij})$
is an alternating $4\times 4$ matrix, this map is given, in a suitable basis of $S^+$, by 
(see Theorem \ref{half-spinor}.\ref{defgamma}):
$$
\gamma:\,Z_B\,\longmapsto\,(z_1:\ldots:z_7)\,=\,
\big(1:b_{12}:b_{13}:b_{14}:b_{12}b_{34}-b_{13}b_{24}+b_{14}b_{23}:
-b_{34}:b_{24}:-b_{23}\big)~.
$$
The coordinate functions are, up to signs, the Pfaffians of the alternating submatrices of $B$
with an even number of rows and columns.
The closure of the image of $\gamma$ is the spinor variety, a smooth quadric:
$$
Q^+\,=\,\gamma(IG(4,V_\CC)^+)\,=\,\{\,(z_1:\ldots:z_8)\,\in\,\PP S^+_\CC:\quad 
z_1z_5+z_2z_6+z_3z_7+z_4z_8\,=\,0\,\}~.
$$
In fact the signs and the order of the coordinate functions on $S^+_\CC$ were chosen in such a 
way as to obtain this simple equation.

The homogeneous coordinates above define a $\ZZ$-module $S^+\cong\ZZ^8\subset S^+_\CC$ 
with bilinear form $(\bullet,\bullet)_{S^+}$ such that for $z=(z_1,\ldots,z_8)\in S^+$ one has
$(z,z)_{S^+}=2(z_1z_5+z_2z_6+z_3z_7+z_4z_8)$.
In particular, $S^+\cong U^4$ and for $z\in S^+_\CC$ one has 
$z\in Q^+$ iff $(z,z)_{S^+}=0$ where we use the $\CC$-bilinear extension of the bilinear form.

\

\subsection{Orthogonal complex structures and their period space $\Omega$ } 
An orthogonal complex structure $J$ on $V_\RR$ is determined by (and determines)
a maximally isotropic subspace $Z_+$ such that $V_\CC=Z_+\oplus\overline{Z_+}$.
Using the spinor map we see that $\ell:=\gamma([Z_+])$ is a point of the quadric
$Q^+\subset\PP S^+_\CC$, that is $(\ell,\ell)_{S^+}=0$. 
Since the spinor map is defined over $\QQ$, we get
$[\overline{Z_+}]=\bar{\ell}$, the complex conjugate of the point $\ell$ in $\PP S^+_\CC$.
The condition that $Z_+\cap \overline{Z_+}=0$ is equivalent to the fact that the complex
line spanned by $\ell,\bar{\ell}$ is not contained in $Q^+$ (see Lemma \ref{acth}).
This again is equivalent to $(\ell,\bar{\ell})_{S^+}\neq 0$ and since
$(\ell,\bar{\ell})_{S^+}\in\RR$ we see that $(\ell,\bar{\ell})_{S^+}$ is either positive or negative.

We define an open (six dimensional, connected) analytic subset of $Q^+$ by
$$
\Omega\,=\,\Omega_{S^+}\,:=\,\{\,\ell\,\in\,\PP S^+_\CC:\;(\ell,\ell)_{S^+}\,=\,0,\quad
(\ell,\overline{\ell})_{S^+}\,>\,0\,\}~.
$$
Then any $\ell\in \Omega$ defines a maximal isotropic subspace $Z_\ell$ of
$V_\CC$ such that $V_\CC=Z_\ell\oplus\overline{Z_\ell}$ 
and thus it defines an orthogonal complex structure $J_\ell$ on $V_\RR$.

The complex structure $J_\ell$ on $V_\RR$ defines a complex torus $\cT_\ell$ of dimension four
by requiring an isomorphism of weight 1 Hodge structures
$$
H^1(\cT_\ell,\ZZ)\,=\,(V,J_\ell),\qquad\mbox{i.e.}\quad
H^{1,0}(\cT_\ell)\,=\,Z_\ell~.
$$
This complex torus can also be defined as $\cT_\ell=V_\CC/(Z_\ell +V)$.

\

\section{Tori with an orthogonal structure and a Cayley class}\label{toricayley}

\subsection{} 
Using representation theory (explained in more detail in the Appendix), we 
recall the relation between the spinor and the Pl\"ucker map. 
We also find a natural map from $S^+$ to $\wedge^4 V$, 
the image of $s\in S^+$ is denoted by $c_s\in \wedge^4V$.
For $\ell\in\Omega$ the complex torus $\cT_\ell$ has $H^1(\cT_\ell,\ZZ)=(V,J_\ell)$. 
Thus we can also identify the Hodge structures 
$\wedge^4V=H^4(\cT_\ell,\ZZ)$ and for $s\in S^+$ 
we obtain a cohomology class $c_s\in H^4(\cT_\ell,\ZZ)$
which is Markman's Cayley class of $s$.

In \S \ref{cs22} we recall Markman's result that the Cayley class is a Hodge class, 
so $c_s\in H^{2,2}(\cT_\ell,\ZZ)$, if and only if 
$\ell\in \Omega_{s^\perp}:=s^\perp\cap \Omega$ where $s^\perp$ is the hyperplane in $S^+_\CC$ 
defined by $s$ using the bilinear form on $S^+$. 
Hence the five dimensional complex manifold $\Omega_{s^\perp}$ parametrizes
the four dimensional complex tori with an orthogonal structure and Hodge class $c_s$.

\subsection{The spinor and the Pl\"ucker map} \label{spinorpluecker}
From the isomorphism $\pi^*\cO_{\PP^N}\cong \cL^{\otimes 2}$ over
$IG(4,V_\CC)^+$, 
one can deduce that the Pl\"ucker map on $IG(4,V_\CC)^+$ is
the composition of the spinor map $\gamma$ with the second Veronese map 
$\nu$ on $\PP S^+_\CC$. The Veronese map is induced by
$$
\nu:\,S^+\,\longrightarrow\,Sym^2(S^+),\qquad s\,\longmapsto\,s\odot s~.
$$
More precisely, the group $Spin(V_\CC)$, a double cover of $SO(V_\CC)$, has a natural 
(half-spin) representation on $S^+_\CC$ and on the 36 dimensional vector space $Sym^2(S^+_\CC)$.
This latter representation is reducible, due to the $Spin(V)$-invariant quadric on $S^+$
which dually defines an invariant one dimensional subspace $\Gamma_0$ of $Sym^2(S_\CC^+)$. 
A complement of this subspace
turns out to be an irreducible $Spin(V_\CC)$-representation and is denoted by $\Gamma_{2\alpha}$:
$$
Sym^2(S^+_\CC)\,\cong\,\Gamma_{2\alpha} \,\oplus\,\Gamma_0~.       
$$
The subspace $\Gamma_{2\alpha}$ is spanned by the symmetric tensors $z\odot z\in Sym^2(S^+_\CC)$
with $[z]\in Q^+\subset \PP S^+_\CC$.

There is a decomposition 
of the $70$-dimensional $\wedge^4V_\CC$ in two
irreducible $Spin(V_\CC)$-representations of dimension $35$ (it corresponds to the decomposition
of $\wedge^4V_\CC$ into dual and anti-selfdual $4$-forms for the Hodge star operator defined by
$(\bullet,\bullet)_V$):
$$
\wedge^4V_\CC\,=\,\Gamma_{2\alpha}\,\oplus\,\Gamma_{2\beta}~.
$$
The image of $Q^+$ by the second Veronese map spans the linear subspace
$\PP \Gamma_{2\alpha}\subset \PP^N=\PP \wedge^4V_\CC$. 

\

\subsection{The Cayley classes}\label{cayleyclasses}
Another consequence of the relation between the $Spin(V)$-representations $Sym^2(S^+)$
and $\wedge^4 V$ 
is that any
element $s\in S^+$ defines a 4-form $c_s\in \wedge^4V$, which is called the Cayley class of $s$
(\cite[Remark 12.4]{Markman}, \cite[\S 2.1]{Munoz}).
It is obtained as the composition
$$
S^+\,\stackrel{\nu}{\longrightarrow}\,Sym^2(S^+)\,\cong\, \Gamma_{2\alpha}\,\oplus\,\Gamma_0
\longrightarrow\,\Gamma_{2\alpha}\,\longrightarrow\,\wedge^4 V~,\qquad s\longmapsto c_s~.
$$
This map is equivariant for the action of $Spin(V)$. The stabilizers in 
$Spin(V)$ of $s$ and $c_s$ thus have the same Lie subalgebra. If $(s,s)_{S_+}\neq 0$
the complexification of this Lie algebra is isomorphic to $so(7)_\CC$.

\subsection{The Cayley class and Hodge classes}\label{cs22}
Let $\ell\in \Omega\subset Q^+$ and let $\cT_\ell$ be the associated complex torus.
The Hodge decomposition of the first cohomology group $H^1(\cT_\ell,\ZZ)=(V,J_\ell)$ is  given
by the eigenspaces $Z_\ell,\overline{Z}_\ell=Z_{\overline{\ell}}$
of the orthogonal complex structure $J_\ell$ in $V_\CC$:
$$
H^1(\cT_\ell,\CC)\,=\,V_\CC\,=\,Z_\ell\,\oplus\,Z_{\overline{\ell}},\qquad
J_\ell\,=\,(i,-i)\,\in\,End(Z_l)\oplus End(Z_{\overline{\ell}})~.
$$
To describe the Hodge structure on $H^k(\cT_\ell,\ZZ)$ we use the homomorphism
$$
h_{V,\ell}:\;U(1)\,:=\,\{z\in\CC:z\overline{z}=1\}\,\longrightarrow\, GL(V_\RR),
\qquad 
h_{V,\ell}(a+bi)\,:=\,aI+bJ_\ell,
$$
where $a,b\in\RR$, $a^2+b^2=1$.
Notice that $aI+bJ_\ell=(a+bi,a-bi)\in End(Z_l)\oplus End(Z_{\overline{\ell}})$.

Since $H^k(\cT_\ell,\ZZ)=\wedge^kH^1(\cT_\ell,\ZZ)=\wedge^kV$, the Hodge decomposition 
$H^k(\cT_\ell,\CC)=\oplus H^{p,q}(\cT_\ell)$ is defined by
$$
H^{p,q}(\cT_\ell)\,=\,\big(\wedge^pZ_\ell\big)\otimes\big(\wedge^qZ_{\overline{\ell}}\big)\,=\,
\{x\,\in\,\wedge^k V_\CC:\;h_{V,\ell}(a+bi)\cdot x\,=\,(a+bi)^p(a-bi)^qx\quad\forall a+bi\in U(1)\,\}~,
$$
In particular, the Hodge classes in $H^{2p}(\cT_\ell,\ZZ)$ are the invariants of the one-parameter
subgroup $h_{V,\ell}$ of $SO(V_\RR)$.
The following proposition 
is essentially \cite[Lemma 12.2]{Markman}.

\subsection{Proposition}\label{CayleyHodge}
Let $c_s\in\wedge^4V$ be the Cayley class defined by $s\in S^+$, the integral lattice, and let
$\ell\in \Omega_{S^+}$. Then $c_s$ is an integral Hodge class in $H^4(\cT_\ell,\ZZ)$ exactly when 
$(\ell,s)_{S^+}=0$:
$$
c_s\in\,H^{2,2}(\cT_\ell,\ZZ)\,:=\,H^4(\cT_\ell,\ZZ)\cap H^{2,2}(\cT_\ell)
\qquad
\mbox{if and only if}\qquad 
\ell\,\in\,\Omega_{s^\perp}\,:=\,\{\ell\in\Omega:\,(\ell,s)_{S^+}\,=\,0\,\}~.
$$

\ts
First we observe that $h_{V,\ell}(z)\in SO(V_\RR)$ for all $z\in U(1)$.
In fact, for $v,w\in V_\RR$ we have
$$
\big((aI+bJ_\ell)v,(aI+bJ_\ell)w\big)_V\,=\,
a^2\big(v,w\big)_V\,+\,ab\Big(\big(v,J_\ell w\big)_V\,+\,\big(J_\ell v,w\big)_V\Big)\,+\,
b^2\big(J_\ell v,J_\ell w\big)_V\,=\,\big(v,w\big)_V~,
$$
because $(J_\ell v,J_\ell w)_V=(v,w)_V$ implies $(v,J_\ell w)=(J_\ell v,J_\ell^2 w)_V$ 
and $J_\ell^2=-I$. 

The homomorphism lifting the one-parameter subgroup 
$h_{V,\ell}:U(1)\rightarrow SO(V_\CC)$ to $Spin(V_\CC)$
is denoted by
$$
h_\ell:\,U(1):=\{z\in\CC:z\overline{z}=1\}\,\longrightarrow\,Spin(V_\CC)~.
$$
The action of $h_\ell(z)\in Spin(V_\CC)$ in the half-spin representation $\rho^+$ on  $S^+_\CC$ is
(see Lemma \ref{acth}):
$$
\rho^+(h_\ell(z))\ell\,=\, z^2\ell,\quad 
\rho^+(h_\ell(z))\overline{\ell}\,=\,\bar{z}^2\overline{\ell},
\qquad 
\rho^+(h_\ell(z))s\,=\,s,\quad\forall s\in\langle \ell,\overline{\ell}\rangle^\perp~.
$$
Using the induced action of  $Spin(V_\CC)$ on $s\odot s\in Sym^2(S^+_\CC)$ and its image 
$c_s\in\wedge^4V_\CC=H^4(\cT_\ell,\CC)$ 
we see that $c_s$ is invariant under $h_\ell(z)$ for all $z\in U(1)$  
if and only if $s$ is invariant, so $s\in \langle \ell,\overline{\ell}\rangle^\perp$.
For $s\in S^+$ the condition $(s,\ell)_{S^+}=0$ implies, by
complex conjugation, that also $(s,\overline{\ell})_{S^+}=0$, which proves the proposition.
\qed

\

\section{Abelian varieties of Weil type}

\subsection{The complex tori $\cT_\ell$ and abelian varieties}
For a point $\ell\in \Omega$, an open subset of the spinor quadric $Q^+$, 
we defined a complex torus
$\cT_\ell$ of dimension four whose first cohomology group is identified with $V$ and 
whose Hodge structure is determined by
$H^{1,0}(\cT_\ell)=Z_\ell$, the maximal isotropic subspace of $V_\CC$ corresponding to $\ell$.

Fixing an $s\in S^+$ we also found that for $\ell \in\Omega_{s^\perp}$ 
this complex torus has an integral Hodge class (the Cayley class) $c_s\in H^{2,2}(\cT_\ell,\ZZ)$.
Now we assume that $(s,s)_{S^+}>0$ and we fix another, non-isotropic, class $h\in s^\perp$ 
with $(h,h)_{S^+}>0$. Hence the rank two sublattice
$\langle h,s\rangle\subset S^+$
generated by $h,s$ is positive definite for the bilinear form on $S^+$. 
For $\ell \in \langle h,s \rangle^\perp\cap \Omega$
the torus $T_\ell$ turns out to be an abelian variety of Weil type and the Cayley class $c_s$
is a non-trivial Hodge class. This result, Theorem \ref{Weil} below, is due to
O'Grady \cite[Theorem 5.1]{O'G} and Markman \cite[Corollary 12.9, Theorem 13.4]{Markman}. 
First we recall the basic facts on abelian varieties of Weil type.

\subsection{Abelian varieties of Weil type}
Let $A$ be an abelian variety and let $K=\QQ(\sqrt{-d})$, with $d\in\ZZ_{>0}$, 
be an imaginary quadratic field. An abelian variety of Weil type (with field
$K$) is a pair $(A,K)$, where $A$ is an abelian variety and $K\hookrightarrow End(A)_\QQ$ 
is a subalgebra of the endomorphism algebra of $A$, such that for all $x\in K$, $x\notin \QQ$,
the endomorphism of $T_0A$
defined by the differential of $x=a+b\sqrt{-d}\in K$, with $a,b\in\QQ$,
has eigenvalues $x=a+b\sqrt{-d}$ and $\bar{x}=a-b\sqrt{-d}$ with the same multiplicity. 
Equivalently, the eigenvalues of $x^*$ on $H^{1,0}$ have the same multiplicity.
In particular, if $(A,K)$ is of Weil type, then $\dim A$ is even.

Given an abelian variety of Weil type $(A,K)$, 
there exists a polarization $\om_K\in H^{1,1}(A,\ZZ)$ on $A$
such that for all $x\in K$ its pull-back is  
$$
x^*\om_K\,=\,\Nm(x)\om_K,\qquad \Nm(x)\,=\,x\bar{x}~,
$$
where $\Nm(x)$ is the norm of $x\in K$ (see \cite[Lemma 5.2.1]{hcav}). 
We call such a 2-form a polarization of Weil type and  
$(A,K,\om_K)$ is called a polarized abelian variety of Weil type.

\subsection{The Weil classes.} \label{weilclasses}
For a general abelian variety of Weil type 
$(A,K)$ of dimension $2n$, the spaces of Hodge classes 
$$
B^p(A)\,:=\,H^{p,p}(A,\QQ)\,:=\,H^{2p}(A,\QQ)\cap H^{p,p}(A)
$$
have dimensions  
(\cite{Weil}, see also \cite[Theorem 6.12]{hcav}):
$$
\dim B^p(A)\,=\,1,\quad (p\ne n),\qquad \dim B^n(A)\,=\,3~.
$$
Since $\dim B^1(A)=1$, any $\om\in B^1(A)$, $\om\ne 0$, 
defines (up to sign) a polarization on $A$ which will be of Weil type.

The action of the multiplicative group $K^\times:=K-\{0\}$ on $H^1(A,K):=H^1(A,\QQ)\otimes_\QQ K$ 
has an eigenspace decomposition into two $2n$-dimensional $K$ subspaces
$$
H^1(A,K)\,=\,Z_\kappa\,\oplus\,Z_{\overline{\kappa}},\qquad
x^* (v,w)\,=\,(xv,\overline{x}w)
$$
that are conjugate over $K$. 
Since $A$ is of Weil type, the complexifications of
these eigenspaces both have Hodge numbers $h^{1,0}=h^{0,1}=n$. 
Thus in $H^{2n}(A,K)=\wedge^{2n} H^1(A,K)$ 
there are two 1-dimensional subspaces 
$\wedge^{2n}Z_\kappa$, $\wedge^{2n}Z_{\overline{\kappa}}$ of Hodge type $(n,n)$.
Since they are conjugate, their direct sum is defined over $\QQ$. 
This defines a 2-dimensional subspace of Hodge classes
$$
W_K\,\subset \,H^{n,n}(A,\QQ),
\qquad 
W_K\otimes_\QQ K\,=\,\wedge^{2n}Z_\kappa\,\oplus\,\wedge^{2n}Z_{\overline{\kappa}}~.
$$
(There is also a natural identification of $W_K$ with $\wedge^{2n}_K H^1(A,\QQ)$ where $H^1(A,\QQ)$
is viewed as a $2n$-dimensional $K$ vector space.) The subspace $W_K$
is called
the space of Weil classes. For any $A$ of Weil type one has
$$
\QQ \omega_K^n\;\oplus\;W_K\;\subseteq\;B^n(A)
$$
where $\om_K^n$, is the $n$-fold exterior product of $\omega_K$ with itself.
For a general $A$ of Weil type one has $B^n(A)=\QQ \omega_K^n\,\oplus\,W_K$.

An element $x\in K$ acts with eigenvalues $(x\bar{x})^n$, $x^{2n}$, $\bar{x}^{2n}$ on
$\QQ \omega_K^n\,\oplus\,W_K$. Thus if a non-zero element $c$ in the three dimensional
$\QQ$ vector space $\QQ \omega_K^n\,\oplus\,W_K$ is algebraic
and it is not an eigenvector for the $K$-action (so it is not a multiple of $\omega^n_K$)
then all classes in $\QQ \omega_K^n\,\oplus\,W_K$ are algebraic since $\omega_K^n$ is and so is
$x^*c$ for all $x\in K$.

\subsection{The Hermitian form}
The $\QQ$ vector space $H_1(A,\QQ)$ is also a $K$ vector space for the action of
$K$ given by $x_*$ for $x\in K\subset End(A)_\QQ$. 
A polarization of Weil type $\omega_K\in H^2(A,\QQ)$ defines an alternating form 
on $H_1(A,\QQ)$ and it also defines
a $K$-valued Hermitian form $H$ on this $K$-vector space by:
$$
H\,:\,H_1(A,\QQ)\times H_1(A,\QQ)\,\longrightarrow\,K,\qquad
H(x,y)\,:=\,\omega_K(x,(\sqrt{-d})_*y)\,+\,\sqrt{-d}\omega_K(x,y).
$$
If $\Psi\in M_n(K)$ is the Hermitian matrix defining $H$ w.r.t.\ some 
$K$-basis of $H_1(A,\QQ)$ 
then $\det(\Psi)\in \QQ^\times=\QQ-\{0\}$ and the class
of $\det(\Psi)\in \QQ^\times/\Nm(K^\times)$, called the discriminant of $H$, 
is independent of the choice of the basis.
Given two non-degenerate
Hermitian forms $H,H'$ on $K^n$, there is a $K$-linear map
$M:K^n\rightarrow K^n$ such that $H(x,y)=H(Mx,My)$ for all $x,y\in K^n$ if and only
if $H$, $H'$ have the same signature and the same discriminant.

The discriminant of a polarized abelian variety of Weil type $(A,K,\omega_K)$ 
is the discriminant of $H$.

In Markman's approach, the real part of $H$, 
which is a bilinear form, is (up to the duality
between $H_1(A,\ZZ)$ and $H^1(A,\ZZ)$ and up to a scalar multiple) 
the bilinear form $(\cdot,\cdot)_V$, cf.\  \S \ref{Thm_pol}. 
In particular, it is the same for all families of Weil type, for all fields, 
considered in \cite{Markman} and in Theorem \ref{Weil} below.

\subsection{Complete families}
Given a $K$ vector space $U$ of dimension $2n$ and a Hermitian form $H:U\times U
\rightarrow K$, any $2n$-dimensional abelian variety of Weil type $A$ with field $K$ 
and discriminant
equal to the discriminant of $H$ is obtained by choosing a free $\ZZ$-module 
$\Lambda\subset U$ of rank $4n$ and a complex structure on 
$J$ on $\Lambda_\RR:=\Lambda\otimes_\ZZ\RR$
such that $J$ commutes with $K$, 
the two eigenspaces of $x\in K$, $x\notin\QQ$, on $(\Lambda_\RR,J)$ have the same
dimension and finally
the imaginary part $\omega_K$ of $H$ 
defines a polarization on the complex torus $(\Lambda_\RR,J)/\Lambda$.

The unitary group $U(H)_\RR\cong U(n,n)$ of the Hermitian form $H$ on the 
$\CC=K\otimes_\QQ\RR$ vector space $\Lambda_\RR$ acts by conjugation $g\cdot J:=gJg^{-1}$
on these complex structures.
From this one obtains a complete family of abelian $2n$-folds of Weil type
parametrized by a Hermitian symmetric domain isomorphic to 
$U(n,n)/(U(n)\times U(n))$, so of complex dimension $n^2$.
The unitary group $SU(H)\subset GL(\Lambda_\QQ)$, viewed as algebraic group over $\QQ$, 
is the special Mumford Tate group of the general
abelian variety in the family, see \cite{hcav}.

\

We discuss the proof of the following theorem
in the remainder of this section.

\subsection{Theorem} \label{Weil}
Let $h,s\in S^+$ be perpendicular and such that $\langle h,s\rangle
\subset S^+$ is a positive definite rank two sublattice. Let $d:=(h,h)_{S^+}(s,s)_{S^+}\in \QQ_{>0}$
and let $\ell \in \Omega_{\{h,s\}^\perp}$, where
$$
\Omega_{\{h,s\}^\perp}\,:=\,\{\ell\in\Omega_{s^\perp}\,:\quad (\ell,h)_{S^+}\,=\,0\,\}\;=\;
\{\ell\in\Omega\,:\quad (\ell,s)_{S^+}\,=\,(\ell,h)_{S^+}\,=\,0\,\}
$$
is a complex manifold of dimension four.
Then we have:

\begin{enumerate}
\item[a)]
The complex four dimensional torus $\cT_\ell$ has endomorphisms by $K=\QQ(\sqrt{-d})$, that is 
$K\subset \mbox{End}(\cT_\ell)_\QQ$.
\item[b)]
The complex torus $\cT_\ell$ has a polarization $\omega_K\in H^2(\cT_\ell,\ZZ)$
and $(\cT_\ell,K,\omega_K)$ is polarized abelian fourfold of Weil type.
\item[c)]
The discriminant of the polarization $\omega_K\in H^2(\cT_\ell,\ZZ)$ is trivial. 
\item[d)]
The Cayley class $c_s\in H^{2,2}(\cT_\ell,\ZZ)$ 
is not contained in the subspace $\QQ\omega_K^2$ 
where $\omega_K^2=\omega_K\wedge\omega_K$.
\item[e)]
The four dimensional family of these fourfolds of Weil type parametrized by 
$\Omega_{\{h,s\}^\perp}$ is complete.
\end{enumerate}

\subsection{Endomorphisms of $\cT_\ell$}\label{endomorphisms}
Since the sublattice $\langle h,s\rangle$ is positive definite, 
we may assume that the restriction $q$ of the quadratic form on $S^+$ 
is given by $q(xh+ys)=ax^2+by^2$, with both $a=(h,h)_{S^+},b=(s,s)_{S^+}\in\QQ$ positive. 
Hence $d=ab>0$. The zero locus of $q$ is defined by $a^{-1}((ax)^2+aby^2)=0$,
showing that there are 
two isotropic lines in $\langle h,s\rangle_\CC$ defined by $ax\pm\sqrt{-d}y=0$.
These two lines are conjugate over $K$ where the conjugation on $K$ is 
$\overline{x+y\sqrt{-d}}=x-y\sqrt{-d}$ with $x,y\in\QQ$.
In $\PP S^+_\CC$ they correspond to the two points of intersection of the line
$\PP\langle h,s\rangle_\CC$ with the spinor quadric $Q^+\cong IG(4,V_\CC)^+$,
which  we denote by $\kappa,\overline{\kappa}$:
$$
\{\kappa,\,\overline{\kappa}\}\,=\,Q^+\,\cap\,\PP\langle h,s\rangle_\CC\qquad(\subset \PP S^+_\CC)~.
$$

As $Q^+=\gamma(IG(4,V_\CC)^+)$, these two points 
define two maximal isotropic subspaces
in $V_K:=V\otimes_\QQ K$ denoted by $Z_\kappa,Z_{\overline{\kappa}}$.
Since the points $\kappa,\overline{\kappa}$ are conjugate over $K$, 
so are these subspaces: if $w=v+\sqrt{-d}v'\in Z_\kappa$ with $v,v'\in V_\QQ$ 
then $\overline{w}=v-\sqrt{-d}v'\in Z_{\overline{\kappa}}$.

The plane $\langle h,s\rangle_\CC$ is not contained 
in $Q^+$, hence 
these two subspaces have trivial intersection (Lemma \ref{acth}, \cite[III.1.12]{chevalley}):
$$
V_K\,=\,Z_\kappa\,\oplus\,Z_{\overline{\kappa}},\qquad 
\overline{(v_1,v_2)}=(\overline{v_2}, \overline{v_1})\qquad (v_1\in Z_\kappa,\,
v_2\in Z_{\overline{\kappa}})~.
$$

We identify the $\QQ$ vector space $V_\QQ$ with the image of $V_\QQ\hookrightarrow V_K$, it
consists of the points $(v_1,\overline{v_1})$ with $v_1\in Z_\kappa$. 
Now we define an action of $K$ on $V_\QQ\, (\subset V_K)$ by
$$
K\times V_\QQ\,\longrightarrow\,V_\QQ,\qquad 
x\cdot (v_1,\overline{v_1})\,:=\,(xv_1,\bar{x}\overline{v_1})\,=\,
(xv_1,\overline{xv_1})
\qquad(\in\,V_\QQ\,\subset\,Z_\kappa\,\oplus\,Z_{\overline{\kappa}})~,
$$
where $\bar{x}$ is the conjugate of $x\in K$. 

To show that this induces an inclusion $K\subset\mbox{End}(\cT_\ell)_\QQ$,
it suffices to verify that any $x\in K$ commutes with the complex structure $J_\ell$ on $V_\RR$.
Since $\ell\in \Omega_{{h,s}^\perp}$ we have $(\ell,\kappa)_{S^+}=0$ and similarly the scalar
products of any one of $\ell,\overline{\ell}$ and any one of $\kappa,\overline{\kappa}$
are zero. Therefore the intersection of $Z_\ell,Z_{\overline{\ell}}$ with the complexifications of
$Z_\kappa,Z_{\overline{\kappa}}$
is not zero by Lemma \ref{acth}. Since these spaces are parametrized by the same connected
component $IG(4,V_\CC)^+$, their intersection is even dimensional and thus it is two dimensional. 
From the eigenspace decomposition for $J_\ell$,
$V_\CC=Z_\ell\oplus Z_{\overline{\ell}}$, we obtain the decomposition
$$
V_\CC\,=\,(Z_\ell\cap Z_{\kappa,\CC})\;\oplus\;(Z_\ell\cap Z_{\overline{\kappa},\CC})\;
\oplus\; (Z_{\overline{\ell}}\cap Z_{\kappa,\CC})\;\oplus\;
(Z_{\overline{\ell}}\cap Z_{\overline{\kappa},\CC})~.
$$
The action of $J_\ell$ and $x\in K$ on these four summands are scalar multiplications (by
$\pm i$ and $x,\bar{x}$ respectively), hence the action of $K$ indeed commutes with $J_\ell$.
Since each summand has dimension $2$, the eigenvalues of $x\in K$, $x\not\in\QQ$, on 
$Z_\ell=H^{1,0}(\cT_\ell)$ have the same dimension.

\subsection{The polarization}\label{Thm_pol}
The combination of the $K$-action on $V_\QQ=H^1(\cT_\ell,\QQ)$  
with the bilinear form $(\bullet,\bullet)_V$ leads a polarization $\omega_K\in H^2(\cT_\ell,\QQ)$ 
on $\cT_\ell$.
We define a bilinear form $E$ on $V_\QQ$ by:
$$
E:\,V\times V\,\longrightarrow\,\QQ,\qquad E(v,w)\,=\,(\sqrt{-d}\cdot v,w)_V~.
$$
The duality $V=H_1(\cT_\ell,\QQ)^{dual}$ implies that $E$ defines an element
$\omega_K\in \wedge^2V=H^2(\cT_\ell,\QQ)$. 
Similar to the computations for K\"ahler forms and metrics we establish the basic properties of $E$
which imply that $(\cT_\ell,K,\omega_K)$ is a polarized abelian fourfold of Weil type.

First of all, we have for all $v,w\in V_\QQ$ and all $x\in K$ that
$$
E(x\cdot v,x\cdot w)\,=\,x\bar{x}E(v,w)~.
$$
To verify this, we extend $E$ $K$-bilinearly to $V_K$ and we use that 
$Z_\kappa,Z_{\overline{\kappa}}$ are isotropic subspaces. Thus, with $v=v_1+\overline{v_1}$,
$w=w_1+\overline{w_1}\in Z_\kappa\oplus Z_{\overline{\kappa}}$ we get
{\renewcommand{\arraystretch}{1.5}
$$
\begin{array}{rcl}
E(x\cdot v,x\cdot w)&=&
\big(x\sqrt{-d}v_1+\overline{x\sqrt{-d}v_1},xw_1+\overline{xw_1})\big)_V\\
&=&
\big(x\sqrt{-d}v_1,\bar{x}\overline{w_1}\big)_V\,+\,\big(\bar{x}\overline{\sqrt{-d}v_1},xw_1\big)_V\\
&=&x\bar{x}\big((\sqrt{-d}v_1,\overline{w_1})_V\,+\,(\overline{\sqrt{-d}v_1},w_1)_V\big)\\
&=&x\bar{x}E(v,w)~.
\end{array}
$$
} 

Next we show that $E$ is alternating:
$$
E(v,w)\,=\,(\sqrt{-d}\cdot v,w)_V\,=\,(w,\sqrt{-d}\cdot v)_V\,=\,
\mbox{$\frac{1}{d}$}(\sqrt{-d}\cdot v,\sqrt{-d}^2\cdot w)_V\,=\,-(\sqrt{-d}\cdot v,w)_V\,=\,-E(w,v)~.
$$

To show that the 2-form $\omega_K$ is of type $(1,1)$ it suffices to show that
$E(J_\ell v,J_\ell w)=E(v,w)$ for all $v,w\in V_\RR$:
$$
E(J_\ell v,J_\ell w)\,=\,\big(\sqrt{-d}\cdot J_\ell v,J_\ell w\big)_V\,=\,
\big( J_\ell(\sqrt{-d}\cdot v),J_\ell w)_V\,=\,\big(\sqrt{-d}\cdot v,w)_V\,=\,E(v,w)~.
$$

Finally we verify that $E(J_\ell v,v)>0$ for non-zero $v\in V_\RR$. That is, we must show that
$(\sqrt{-d}\cdot J_\ell v,w)>0$. The endomorphisms $\sqrt{-d},J_\ell$ of $V_\RR$ are both
constructed from decompositions of $V_\CC$ with two conjugate isotropic subspaces 
$Z_\kappa, Z_{\overline{\kappa}}$ and $Z_\ell,Z_{\overline{\ell}}$ respectively. 
The corresponding points $\kappa,\overline{\kappa},\ell,\overline{\ell}\in Q^+=IG(4,V)^+$
span a $\PP^3\in \PP S^+_\CC$ which is the projectivization of the complexification of 
the four dimensional subspace $<h,s,\ell+\overline{\ell},(\ell-\overline{\ell})/i>\subset
S^+_\RR$ (here $\CC=\RR+i\RR$). Notice that this basis consists of perpendicular vectors
for $(\bullet,\bullet)_{S^+}$ and that the subspace is positive definite.

The group $Spin(V_\RR)$ acts via $SO(S^+_\RR)$ on $S^+_\RR$ and this action is 
transitive on such  subspaces. As $Spin(V_\RR)$ 
also acts via $SO(V_\RR)$ on $V_\RR$, we see that it
suffices to show that $(J_1J_2v,v)>0$ for all non-zero $v\in V_\RR$ where the linear maps
$J_1,J_2$ are defined by any two orthogonal positive definite 2-dimensional subspaces of $S^+_\RR$.
(Markman shows that the map $J_1J_2$ is already determined, up to a scalar multiple, 
by the direct sum of these subspaces.)

We use the conventions from \S \ref{spinormap}. A point 
$z=(z_1,\ldots,z_8)\in S^+\cong U^4$  will be written as
$$
z\,=\,\left(\begin{pmatrix}z_1\\z_5\end{pmatrix},\begin{pmatrix}z_2\\z_6\end{pmatrix},
\begin{pmatrix}z_3\\z_7\end{pmatrix},\begin{pmatrix}z_4\\z_8\end{pmatrix}\right),\qquad
(z,z)_{S^+}=2(z_1z_5+\ldots+z_4z_8)~.
$$
The following four points $\nu_1,\ldots,\nu_4$, where $\nu:=(^1_1)\in U$,
in $S^+$ are perpendicular and span a positive 4-plane in $S^+_\RR$ since $(\nu_i,\nu_i)_{S^+}=8$
and we also define $\ell_1,\ell_2\in S^+_\CC$:
{\renewcommand{\arraystretch}{1.5}
$$
\begin{array}{rcl}
\nu_1&=&(\nu,\phantom{-}\nu,\phantom{-}\nu,\phantom{-}\nu),\\
\nu_2&=&(\nu,\phantom{-}\nu,-\nu,-\nu),\\
\nu_3&=&(\nu,-\nu,\phantom{-}\nu,-\nu),\\
\nu_4&=&(\nu,-\nu,-\nu,\phantom{-}\nu),\\
\end{array}
\qquad
{\renewcommand{\arraystretch}{2.5}
\begin{array}{rccccl}
\ell_1&:=&(\nu_1+i\nu_2)/(1+i)&=&(\nu,\phantom{-}\nu,-i\nu,-i\nu)\\
\ell_2&:=&(\nu_3+i\nu_4)/(1+i)&=&(\nu,-\nu,-i\nu,\phantom{-}i\nu)\\
\end{array}
}
$$
}
Then $\ell_1,\overline{\ell_1}$ and $\ell_2,
\overline{\ell_2}$ are all isotropic vectors and they span
$\langle\nu_1,\nu_2\rangle_\CC$ and $\langle\nu_3,\nu_4\rangle_\CC$ respectively.
Isotropic vectors are in $Q^+=\gamma(IG(4,V_\CC)^+)$ and since these four
all have first coordinate $z_1=1$ they are in the image of the 
open set $IG(4,V_\CC)^+_0$ parametrized by the alternating $4\times 4$ matrices. 
Using the explicit description of $\gamma$ one finds 
$$
\ell_k\,=\,\gamma(Z_{B_k})\quad (k=1,2),\qquad
B_1\,=\,\begin{pmatrix}0&1&-i&-i\\-1&0&i&-i \\i&-i &0&-1\\i&i&1&0\end{pmatrix},\qquad
B_2\,=\,\begin{pmatrix}0&-1&-i&i\\1&0&-i&-i \\i&i &0&1\\-i&i&-1&0\end{pmatrix}~.
$$
The eigenspace with eigenvalue $-1=i^2=(-i)^2$ of the endomorphism $J_1J_2$ of $V_\RR$ 
is the direct sum of $Z_{\ell_1}\cap Z_{l_2}$ and its complex conjugate.
Let $c_k,d_k$ denote the $k$-th column of the matrix $(^{B_1}_{\,I})$, $(^{B_2}_{\,I})$
respectively, then $Z_{\ell_1},
Z_{\ell_2}$ are spanned by the $c_k$ and the $d_k$ ($k=1,\ldots,4$) respectively.
Their intersection is spanned by
$$
c_1-ic_3\,=\,{d_1-id_3},\quad c_2-ic_4\,=\,{d_2-id_4}\qquad
(\in Z_{\ell_1}\cap Z_{\ell_2})~.
$$
Considering $(c_1-ic_3)\pm\overline{(c_1-ic_3})$ etc.,
one finds a basis of the $-1$-eigenspace of $J_1,J_2$. Its perpendicular is the $+1$-eigenspace.
Recall that $e_1,\ldots,e_8$ are the basis vectors of $V$ as in \ref{defV}, then the eigenspace
decomposition is:
$$
V_\RR\,=\,V_+\,\oplus\,V_-\,=\, \langle e_1+e_5,e_2+e_6,e_3+e_7,e_4+e_8\rangle_\RR\,\oplus\,
\langle e_1-e_5,e_2-e_6,e_3-e_7,e_4-e_8\rangle_\RR~.
$$
Notice that $(\bullet,\bullet)_V$ is positive definite on $V_+$ and negative definite on $V_-$.
Writing $v=v_++v_-$ as sum of $J_1J_2$ eigenvectors, one has $(J_1J_2v,v)_V=(v_+,v_+)_V-(v_-,v_-)_V$
and thus indeed $(J_1J_2v,v)_V>0$ for all non-zero $v\in V_\RR$.

\subsection{The discriminant} We refer to \cite[Lemma 12.11]{Markman} (cf.\ \cite[Theorem 5.1]{O'G}) 
for the computation of the
discriminant. See also Proposition \ref{proptrivdisc} for a proof of the triviality of the 
discriminant using results from Lombardo \cite{Lombardo}.

\subsection{The Cayley class and the Weil classes} \label{ccwc}
We define two subgroups of $Spin(V)$.
Let $Spin(V)_s$ be the subgroup which fixes $s\in S^+$ 
and let $Spin(V)_{h,s}$ be the subgroup
which fixes all elements in $\langle h,s\rangle$.
Then one can show that the Cayley class
$c_s$ is the unique $Spin(V)_s$-invariant in $\wedge^4 V$ and that $\omega_K$ is the unique
$Spin(V)_{s,h}$-invariant in $\wedge^2V$. This implies that $c_s\not\in\QQ \omega_K^2$
(cf.\ \cite[Prop 2]{Munoz}, \cite[Thm 13.4]{Markman} and \S \ref{spin7}). 

One can also use that the $K\otimes_\QQ\CC \cong\CC\times\CC$-action 
on $V_\CC$ has the eigenspaces $(Z_\kappa)_\CC,(Z_{\bar{\kappa}})_\CC$. The one parameter
subgroup $h_R$ of $Spin(V_\CC)$ which acts as multiplication by $t,t^{-1}$ respectively 
on these eigenspaces fixes $E$, and thus it fixes $\omega_K\in \wedge^2 V$ and also 
$\omega^2_K\in \wedge^4V$. 
On the other hand, $h_R$ has eigenvalues $t^2,t^{-2}$ on 
$\langle \kappa,\bar{\kappa}\rangle_\CC=\langle h,s\rangle_\CC\subset S^+_\CC$ 
by Lemma \ref{acth}. Therefore $c_s$, the image of $s\odot s$ in $\wedge^4V$, 
is not invariant under $h_R$ and thus it
cannot be a multiple of $\omega^2_K$.

\subsection{Complete families}
The Lie group $Spin(V_\RR)_{h,s}$ acts  on $\Omega_{\{h,s\}^\perp}$.
This action induces an action of $Spin(V_\RR)_{h,s}$ on the orthogonal complex 
structures on $V_\RR$ by $J_{g\cdot \ell}=gJ_\ell g^{-1}$. 
The fixed points $\kappa,\bar{\kappa}\in Q^+\cap \langle h,s\rangle_\CC$ of the action
of $Spin(V_\RR)_{h,s}$ on $Q^+$ correspond to the eigenspaces 
$Z_{\kappa,\CC},Z_{\bar{\kappa},\CC}$ of the $K$-action,
which are thus mapped into themselves. 
This implies that the image of $Spin(V_\RR)_{h,s}$ in $SO(V_\RR)$ commutes 
with the $K$ action on $V_\RR$. 
This image thus preserves the Hermitian form $H$ and therefore $Spin(V)_{h,s}$ maps to the 
algebraic group $SU(H)$ which is the Mumford Tate group of the general $\cT_\ell$ with
$\ell\in \Omega_{\{h,s\}^\perp}$.
For dimension reasons this map is surjective on the real points of these groups and thus
the family of abelian fourfolds of Weil type is complete.

 \

\section{Moduli spaces of sheaves on an abelian surface}\label{moduli}

\subsection{} The constructions considered thus far have a natural geometrical interpretation
in terms of moduli spaces of sheaves on abelian surfaces. We now briefly recall the basic 
definitions and results, due to Mukai and Yoshioka. 
The notation used thus far is now adapted to this context, for example, the free $\ZZ$-module
$W$ of rank four will become $W=H^1(X,\ZZ)$ for an abelian surface $X$ etc.

We conclude with a 
brief outline of Markman's proof of the Hodge conjecture for the general abelian
fourfolds of Weil type with trivial discriminant.

\subsection{The Mukai lattice of an abelian surface}\label{MukaiLattice}
Let $X$ be an abelian surface and let $\hat{X}=Pic^0(X)$ be the dual abelian surface. Let 
$$
W\,=\,H^1(X,\ZZ),\quad W^*\,=\,H^1(\hat{X},\ZZ)\,=\,H^1(X,\ZZ)^*,\qquad V\,:=\,W\oplus W^*~.
$$
The Chern character of a coherent sheaf on $X$ takes values in 
$$
S^+\,:=\,\wedge^{even} H^*(X,\ZZ)\,=\,H^0(X,\ZZ)\,\oplus\,H^2(X,\ZZ)\,\oplus\,H^4(X,\ZZ)~,
$$
and we will identify $H^0(X,\ZZ)$, $H^4(X,\ZZ)$ with $\ZZ$, using the generators $1$ and a 
volume form compatible with the orientation on the complex manifold $X$.

The Mukai lattice of $X$ is the (free, rank $8$) $\ZZ$-module $S^+$ with the bilinear form
given by (this bilinear form coincides up to sign with $(\bullet,\bullet)_{S^+}$):
$$
(r,c,s)\cdot (r',c',s')\,:=\,-(rs'+r's)+\,c\wedge c'~.
$$

For $v=(r,c,s)\in S^+$, with $r>0$, $c\in NS(X)\subset H^2(X,\ZZ)$ and $v^2\geq 6$ the moduli space
of sheaves $E$ on $X$ with $ch(E)=v$, denoted by $\cM(v)$,
is a smooth holomorphic symplectic manifold of dimension $v^2+2$.

\subsection{The case $v=s_n$}\label{vsn}
We now take $v=s_n=(1,0,-n)$, so that $v^2=2n\geq 6$ and $\dim M(v)=2n+2$.
Let $Z\subset X$ be a subscheme of length $n$, then its ideal sheaf $\cI_Z$ has $ch(\cI_Z)=v$
(for an abelian surface, the Chern character $ch(E)$ is the Mukai vector $v(E)$ of the sheaf $E$).
This induces an inclusion of complex manifolds 
$$
Hilb^n(X)\,=\,X^{[n]}\;\hookrightarrow\; \cM(v)\qquad (v=s_n=(1,0,-n))~.
$$
For $\cL\in \hat{X}$ and $\cI_Z\in X^{[n]}$ one also has $\cL\otimes \cI_Z\in \cM(v)$.

The Albanese map $\alpha:X^{[n]}\rightarrow X$ of $X^{[n]}$ fits in a diagram:
$$
\begin{array}{rcl}
 X^{[n]}&&\\
 \downarrow&\searrow \alpha&\\
 X^{(n)}&\stackrel{\Sigma}{\longrightarrow}&X
\end{array}\qquad
\Sigma([p_1,\ldots,p_n])\,:=\,p_1+\ldots+p_n~,
$$
here $X^{(n)}$ is the $n$-th symmetric power of $X$ and $[p_1,\ldots,p_n]\in X^{(n)}$
is the image of $(p_1,\ldots,p_n)\in X^n$ in $X^{(n)}$.

The generalized Kummer variety $K_{n-1}(X)$, of dimension $2n-2$, 
is the irreducible 
holomorphic symplectic manifold obtained as 
$$
K_{n-1}(X)\,=\,\alpha^{-1}(0)\subset X^{[n]}~.
$$

Using locally free resolutions of sheaves one defines a determinant map 
$\det:M(v)\rightarrow \hat{X}$ and one has $\det(\cL\otimes \cI_Z)=\cL$ for $\cL\in \hat{X}$. 
Yoshioka \cite{Yoshioka} showed that
$$
M(v)\,\cong\,\hat{X}\times ({\det}^{-1})(\cO_X)\,\cong\,\hat{X}\,\times\,X^{[n]}\,\cong\,
\hat{X}\times\Big((X\times K_{n-1}(X))/X[n]\Big)
$$
where $X[n]\subset X$ is the subgroup of $n$-torsion points. In particular, the Bogomolov
decomposition of $M(v)$ is the product of the abelian fourfold $X\times\hat{X}$ and the 
irreducible holomorphic symplectic manifold $K_{n-1}(X)$.

\subsection{The cohomology of the generalized Kummer variety}
The composition of the Mukai homomorphism and the restriction map 
$$
v^\perp\,\longrightarrow\,H^2(M(v),\ZZ)\,\longrightarrow\,H^2(K_{n-1}(X),\ZZ)
$$
induces a Hodge isometry (for the weight two Hodge structure on $v^\perp$
defined by $(v^\perp)^{2,0}=H^{2,0}(X)$ and
with the BBF quadratic form on $H^2(K_{n-1}(X),\ZZ)$)
\cite[Thm.\ 0.2]{Yoshioka}.

This implies, by the surjectivity of the period map and with $v=s_n=s$, 
that $\Omega_{s^\perp}$ is the period space of deformations of $K_{n-1}(X)$,
these deformations are called Kummer type varieties.

Moreover, $h^{3,0}(K_{n-1}(X))=0$ so that $H^3(K_{n-1}(X),\CC)=H^{2,1}\oplus H^{1,2}$ is essentially
the first cohomology group of its intermediate Jacobian $H^3(\CC)/(H^{2,1}\oplus H^3(\ZZ))$
and one has (\cite[Prop.\ 4.20]{Yoshioka}):
$$
H^3(K_{n-1}(X),\ZZ)\,=\,H^1(X,\ZZ)\,\oplus\,
H^3(X,\ZZ)\,\cong\,H^1(X,\ZZ)\,\oplus\,H^1(\hat{X},\ZZ)\,=\,V~.
$$

O'Grady and Markman showed that for $\ell \in \Omega_{s^\perp}$ and any deformation 
$Y_\ell$ of $K_{n-1}(X)$ with period $H^{2,0}(Y_\ell)=\CC\ell\subset
(s^\perp)_\CC$, there is an isomorphism of Hodge structures (up to Tate twist and isogeny)
$H^3(Y_\ell,\ZZ)=H^1(\cT_\ell,\ZZ)$. 
In case the complex manifold $Y_\ell$ is algebraic and $h\in H^2(Y,\ZZ)=s^\perp$
is the class of an ample divisor, hence $\ell\in\Omega_{\{h,s\}^\perp}$, 
O'Grady \cite{O'G} showed that the torus $\cT_\ell$ is an abelian variety of Weil type.
Moreover, he showed that for algebraic $Y_\ell$ the Kuga Satake variety of the weight two polarized
Hodge structure of rank six $h^\perp\subset H^2(Y_\ell,\ZZ)$ is (isogeneous to) $\cT_\ell^4$
(see also \S \ref{H2ell} where $h^\perp\cong H^2_\ell$). 

O'Grady also makes a detailed study of the cohomology of generalized
Kummer varieties and in particular he observes that there is a natural map 
(recall $\dim Y_\ell=\dim K_{n-1}(X)=2n-2$):
$$
H^3(Y_\ell,\ZZ)\,\longrightarrow\, H^{4n-6}(Y_\ell,\ZZ)\,\longrightarrow\, H^2(Y_\ell,\ZZ)^{\vee},
$$
the last map is Poincar\'e duality,
which relates the Hodge structures on $H^3(Y_\ell)$ and $H^2(Y_\ell)$.

\subsection{Markman's theorem} Given a sheaf $F'\in \cM(v)$ ($v=s_n$ as in \S \ref{vsn}),
there is a natural map 
$$
\iota_{F'}\,:\,X\times\hat{X}\,\longrightarrow\, \cM(v),\qquad 
(x,\cL)\,\longmapsto \, (t_x^*F')\otimes \cL
$$
where $t_x:X\rightarrow X$, $y\mapsto x+y$ is the translation by $x$.  
Deforming $K_{n-1}(X)$ to $Y_\ell$,  
with $\ell\in\Omega_{s^\perp}$, 
this map deforms to a map
$$
\iota\,:\,\cT_\ell\,\longrightarrow\,Y_\ell.
$$

A universal sheaf $\cE$ on $X\times \cM(v)$ defines a sheaf $E$ on $M(v)\times M(v)$
by $E:=\mathcal{E}xt^1_{\pi_{13}}(\pi^*_{12}\cE,\pi_{23}^*\cE)$ where $\pi_{ij}$ are the projections
from $M(v)\times X\times M(v)$. For $F\in M(v)$ let $E_F$ the restriction of $E$ to 
$\{F\}\times M(v)=M(v)$. This defines a sheaf on $X\times \hat{X}$ 
whose second Chern class is exactly the Cayley class defined by $v=s_n\in S^+$
(\cite[Prop.\ 11.2]{Markman}, see also \ Prop.\ \ref{explicitCs}):
$$
c_2(\iota_{F'}^*{\mathcal End}(E_F))\,=\,c_v\,\in\,\wedge^4V\,=\,H^4(X\times \hat{X},\ZZ)~.
$$

Markman, using results of Verbitsky,
shows that the sheaf $E_F$ on $M(v)$ deforms to a sheaf 
over any deformation $Y_\ell$ of $M(v)$. Thus $c_v\in H^4(\cT_\ell,\ZZ)$ is an algebraic
class whenever $\cT_\ell$ is an abelian variety. From Theorem \ref{Weil}.d we have that
$c_v$ is not an eigenvector for the action of the multiplicative group $K^\times$ on the Hodge
classes in 
$\QQ\omega^2_K\oplus W_K\subset H^{2,2}(\cT_\ell,\ZZ)$. Thus $\omega_K^2,c_v$ and the images
of $c_v$ under the $K^\times$ action span $\QQ\omega^2_K\oplus W_K$. Since any fourfold 
of Weil type with trivial discriminant is isogeneous to a $T_\ell$, for any
such fourfold the space $W_K$ is spanned by algebraic classes.

\quad

\section{Appendix: The spinor map}\label{App}
\subsection{Background} \label{background}
The spinor map was defined by Cartan \cite{Cartan} (see also \cite{BaloghHH}).
The description given by Chevalley in \cite{chevalley} 
was used by Markman \cite[\S 2]{ManivelSpinor}. We define the spinor map
using the representation theory of orthogonal groups as in \cite[Chapter 20]{FH}
(but our $(v,w)_V$ is $2Q(v,w)$ in \cite{FH}). 

We change the notation: in this Appendix $V$ stands for $V_\CC$, $S^+$ for $S^+_\CC$ etc.
so all $\ZZ$-modules are replaced by their complexifications. 
Whenever convenient we will also write $\CC^{2n}$ for $V$ and $SO(2n)$ for $SO(V)$ etc.

\

\subsection{The Clifford algebra of $V$}\label{cliffordV}
The Clifford algebra $C(V)$ of the complex vector space $V$, of dimension
$2n$, with the bilinear form
$(\bullet,\bullet)_V$ is the quotient of the tensor algebra
$$
C(V)\,:=\,\oplus_{k\geq 0} 
V^{\otimes k}/\langle v\otimes w\,+\,w\otimes v\,-\,(v,w)_V\cdot 1\,\rangle~, 
$$
by the two sided ideal
generated by the $v\otimes w+w\otimes v-(v,w)_V$ with $v,w\in V$, or equivalently,
by the two sided ideal generated by the $v\otimes v-(1/2)(v,v)_V$ for $v\in V$.

The Clifford algebra has dimension $2^{2n}$. We identify $V$ with its image in $C(V)$. 
The even Clifford algebra $C(V)^+$ is the image of $\oplus_{k\geq 0} V^{\otimes 2k}$.

\

Let $V=W\oplus W^*$ be the complexification of the lattices in \S \ref{defV}.
Since $W,W^*$ are isotropic one has $vw=-wv\in C(V)$ for all $v,w\in W$ and 
also for all
$v,w\in W^*$. The subalgebras of $C(V)$ generated by $W,W^*$ are isomorphic to the exterior 
algebras $\wedge^\bullet W$ and $\wedge^\bullet W^*$ respectively. 

Let $e^*:=e_{n+1}\cdots e_{2n}\in C(V)$ be the product of the elements in a basis of $W^*$.
Then the left ideal $S:=C(V)e^*$ of $C(V)$ is isomorphic, as a $\CC$ vector space, 
to $\wedge^\bullet W$, 
$$
\sigma:\,\wedge^\bullet W\,\stackrel{\cong}{\longrightarrow}\,S\,:=\,C(V)e^*,
\qquad w_1\wedge w_2\wedge \ldots \wedge w_r\,\longmapsto\,w_1w_2\ldots w_re^*~,
$$
(\cite[II.2.2]{chevalley}, \cite[Exercise 20.12]{FH}). Under this isomorphism, left multiplication
by $w\in W$ and $w^*\in W^*$ on $S$ correspond to the following endomorphisms of $\wedge^\bullet W$:
$$
w\sigma(\eta)\,=\,\sigma( w\wedge \eta),\qquad 
w^*\sigma(\eta)\,=\,\sigma(D_{w^*}\eta)~,\qquad 
(\eta\in\,\wedge^\bullet W)~,
$$
where $D_{w^*}$ is the derivation on $\wedge^\bullet W$ defined by 
$$
D_{w^*}(1)=0,\qquad  D_{w^*}(w_1\wedge\ldots\wedge w_r)\,=\,\sum_{i=1}^r(-1)^{i-1}
w^*(w_i)(w_1\wedge\ldots\wedge\widehat{w_i}\wedge\ldots\wedge w_r)~,
$$
(here $w^*(w)=(w,w^*)_V$ for $w\in W,w^*\in W^*$).

These operations of $W,W^*$ on $\wedge^\bullet W$ define a $C(V)$-module structure and
$\sigma$ is a homomorphism of $C(V)$-modules. It induces an isomorphism of $\CC$-algebras
between the even Clifford algebra and the direct sum of two matrix algebras (cf.\ \cite[(20.13)]{FH})
$$
C(V)^+\,\cong\,\mbox{End}(S^+)\,\oplus\,\mbox{End}(S^-),\qquad S^+\,:=\,\wedge^{even}W,\qquad
S^-\,:=\,\wedge^{odd}W~.
$$
Since $\dim W=n$ one has $\dim S^\pm=2^{n-1}$.

\subsection{The spin group of $V$} \label{SpinV}
The conjugation on $C(V)$ is the anti-involution given by
$$
x\,:=\,x_1\cdots x_r\,\longmapsto\,x^*\,:=\,(-1)^rx_r\cdots x_1,
$$
notice that it maps $C(V)^+$ into itself. The spin group of $V$ is
$$
Spin(V)\,:=\,\{x\in C(V)^+:\;xx^*\,=\,1,\quad xVx^*\,\subset\,V\,\}~. 
$$
Elements in $Spin(V)$ thus induce linear maps on $V$ and one has the following result.

\subsection{Theorem} There is a surjective homomorphism of complex Lie groups
$$
\rho_V:\,Spin(V)\,\longrightarrow\,SO(V),\qquad x\,\longmapsto [v\,\longmapsto\, xvx^*]
$$
with kernel $\{\pm 1\}$.

\ts For a proof see \cite[Thm 20.28]{FH}. 
\qed

\

\subsection{The half-spin representations}
Besides this `standard representation' of $Spin(V)$ on $V$, 
one also has the two half-spin representations $\rho^+,\rho^-$ of $Spin(V)$ 
on $S^+$ and $S^-$ respectively (vector spaces of dimension $2^{n-1}$), 
given by left multiplication in $C(V)$:
$$
\rho^\pm:\,Spin(V)\,\longrightarrow\,GL(S^\pm),\qquad x\,\longmapsto [\eta\,\longmapsto\, x\eta]~.
$$
See \cite[Exercise 20.38]{FH} for the fact that for $n\equiv 0\mod 4$
the image of $Spin(V)$ lies in $SO(2^{n-1})$
(for a certain bilinear form $\beta$ on $S^+\subset \wedge^\bullet W$ also considered in
\cite[3.2]{chevalley}). The center of $Spin(V)$, $\dim V>2$, 
is isomorphic to $(\ZZ/2\ZZ)^2$ if $n$ is even and is 
cyclic of order four otherwise (cf.\ \cite[Exercise 20.36]{FH}). 
For $n$ even, $n>2$, the three  quotients of $Spin(V)$ by the order two subgroups of
the center are $SO(V)$ and the images of $Spin(V)$ in the two half-spin representations.

\subsection{The Lie algebra $spin(V)= so(2n)$}
The Lie algebra $spin(V)$ of the subgroup $Spin(V)$ of the multiplicative group of $C(V)^+$ 
consists of the $x\in C(V)^+$ such that $x+x^*=0$ and
$xv-vx\in V$ for all $v\in V$ (cf.\ \cite[p.67-68]{chevalley}). It has a basis consisting of
the following $n(n-1)/2+n(n-1)/2+n^2=n(2n-1)$ elements:
$$
e_ie_j,\quad e_{i+n}e_{j+n}\quad  \mbox{with}\quad 1\leq i\leq j\leq n;\qquad 
e_ie_{j+n}\,-\,\mbox{$\frac{1}{2}$}1,\quad 1\leq i,j\leq n~.
$$
To see that these elements are in $spin(V)$ (and to find their action on $V$) one can use
that for $x,y,v\in V$ one has
$$
[xy,v]\,:=\,xyv\,-\,vxy\,=\,x(-vy+(y,v)_V)\,-\,(-xv+(x,v)_V)y\,=\,(y,v)_Vx\,-\,(x,v)_Vy~.
$$

The Lie algebra $spin(V)$ is isomorphic to 
the Lie algebra $so(2n)$ of the orthogonal group $SO(V)=SO(2n)$. This Lie algebra consists of 
the $X\in \mbox{End}(V)$ such that $(Xv,w)_V+(v,Xw)_V=0$ for all $v,w\in V$. 
One finds that  
$$
so(2n)\,=\,\left\{\,X\,=\,\begin{pmatrix} A&B\\C&D\end{pmatrix}\in\mbox{End}(V)\;:\quad 
A\,=\,-{}^tD,\quad {}^tB\,=\,-B,\quad {}^tC\,=\,-C\;\right\}~.
$$
An isomorphism $spin(V)\rightarrow so(2n)$ is given by the differential of $\rho_V$,
so by the representation of $spin(V)$
on $V$ given by $x\cdot v:=xv-vx$. 
Using the computation of $[xy,v]$ above, 
one verifies that this isomorphism is given by
$$
spin(V)\,\stackrel{\cong}{\longrightarrow}\, so(2n),\qquad
\left\{
\begin{array}{rclcl}
e_ie_{n+j}&\longmapsto&X_{i,j},&\qquad&X_{i,j}\,:=\,E_{i,j}\,-\,E_{n+j,n+i}, \\
e_ie_j&\longmapsto&Y_{i,j},&&Y_{i,j}\,:=\,E_{i,n+j}\,-\,E_{j,n+i}, \\
e_{i+n}e_{j+n}&\longmapsto&Z_{i,j},&&Z_{i,j}\,:=\,E_{n+i,j}\,-\,E_{n+j,i}~. \\
\end{array}\right.
$$

We choose the Cartan subalgebra of $so(2n)$ to be the diagonal matrices in $so(2n)$
(as in \cite[\S 18.1]{FH}):
$$
\mathfrak{h} \,:=\,\oplus_{i=1}^n\,\CC H_i,\qquad H_i\,:=\,E_{i,i}\,-\,E_{n+i,n+i}~.
$$
The dual $\mathfrak{h}^*$ of $\mathfrak{h}$ then consists of the linear maps (weights)
$$
\mathfrak{h}^*\,:=\,\oplus_{i=1}^n \CC L_i,\qquad L_i(\sum_{j=1}^n t_jH_j)\,:=\,t_i~.
$$

\subsection{The spinor map}
As before, we identify
$$
V\,=\,W\,\oplus\,W^*\,\qquad 
W\,:=\,\langle e_1,\ldots,e_n\rangle,\quad 
W^*\,:=\,\langle e_{n+1},\ldots,e_{2n}\rangle,
$$
with $W^*=\mbox{Hom}(W,\CC)$ the dual of $W$, where $w^*(w):=Q(w,w^*)$ for $w\in W,w^*\in W^*$.
We denote by $IG(n,2n)^+$ the connected component of the Grassmannian of maximal isotropic subspaces
of $V=\CC^{2n}$ that contains $W^*$. A complex maximally isotropic subspace $Z$ defines a point 
$[Z]\in IG(n,2n)^+$ 
if and only if $\dim(Z\cap W_\CC^*)\equiv n\mod 2$ is even.

We denote by $Z_B$, for an  alternating
$n\times n$ matrix $B$, the maximal isotropic subspace spanned by the columns of $(^B_I)$
analogous to \S \ref{spinormap}, notice that $W^*=Z_0$.

The Grassmannian $IG(n,2n)^+=SO(V)/P$ is a homogeneous space where $P=P_{W^*}$ is
the stabilizer of $W^*$ in the group $SO(2n)$. 
The Lie algebra of $P$, 
which are the $X\in so(2n)$ with $XW^*\subset W^*$, 
consists of the $X\in so(2n)$ with $B=0$.

We recall that the Pfaffian of an alternating $2m\times 2m$ matrix $A$ 
is the complex number $\mbox{Pfaff}(A)$
defined by the following identity in $\wedge^{2m}\CC^{2m}$:
$$
\mbox{Pfaff}(A)e_1\wedge\ldots\wedge e_{2m}
\,=\, m!\omega_A^m,\qquad (\omega_A\,:=\,\sum_{i<j}a_{ij}e_i\wedge e_j)~.
$$

\subsection{Theorem}\label{half-spinor}
Let $\rho^+:Spin(V)\rightarrow GL(S^+)$ be the  half-spin representation of 
$Spin(V)$ on $S^+=\wedge^{even} W:=\oplus_k \wedge^{2k}W$.

\begin{enumerate}
\item In case $n$ is even, the highest weight of $S^+$ is $(L_1+\ldots+L_n)/2$
and it is $(L_1+\ldots+L_{n-1}-L_n)/2$ if $n$ is odd.
\item
The one dimensional subspace
$$
\langle 1\rangle\,=\langle \wedge^0W\rangle\,\subset\,\wedge^{even} W
$$
is invariant under the Lie algebra of $P$. Thus there is a $Spin(V)$ equivariant map
$$
\gamma:\,IG(n,2n)^+\,\longrightarrow\,\PP S^+,\qquad \gamma([\rho_V(\tilde{g})W^*])\,=\,
\rho^+(\tilde{g})1
$$
for $\tilde{g}\in Spin(V)$. 
\item
For an alternating matrix $B\in M_n(\CC)$, let 
$$
X_B\,:=\,\begin{pmatrix}0&B\\0&0\end{pmatrix}\in so(2n),\qquad
\tilde{g}_B\,:=\,exp(X_B)\in Spin(V)~.
$$
In the standard representation $\rho_V:Spin(V)\rightarrow SO(2n)$ one has
$$
\rho_V(\tilde{g}_B)\,=\,\begin{pmatrix}I&B\\0&I\end{pmatrix}\;
(\in\,SO(2n))~\qquad\mbox{and}
\quad 
\rho_V(\tilde{g}_B)Z_0\,=\,Z_B~.
$$
In the half-spin representation on $S^+$ the action of $\tilde{g}_B$
is given by a left multiplication:
$$
\rho^+(\tilde{g}_B):S^+\,\longrightarrow \,S^+,\qquad 
\omega\,\longmapsto \, exp(\omega_B)\wedge\omega~,
$$
and one has
$$
exp(\omega_B)\,=\,\sum_{I,\sharp I\equiv 0\mod 2} \mbox{Pfaff}(B_I)e_I~,
$$
where $I$ runs of over the subsets of $\{1,\ldots,n\}$ with an even number of elements
and $e_I=e_{i_1}\wedge\ldots\wedge e_{i_{2k}}\in\wedge^{even} W=S^+$ with $i_1<\ldots<i_{2k}$.
\item \label{defgamma}
In the basis of $S^+$ consisting of the $e_I$,
the spinor map $\gamma$ on the open subset $IG(n,2n)^+_0$ is given by 
$$
\gamma:\,IG(n,2n)^+_0\,\longrightarrow\,\PP S^+,\qquad [Z_B]\,\longmapsto\,
(\ldots:\mbox{Pfaff}(B_I):\ldots)~.
$$
The image of $\gamma$ is defined by quadrics.
\end{enumerate}

\ts 
The highest weight of the half-spin representation $S^+$ is determined in
\cite[Proposition 20.15]{FH}. 

The Lie algebra of $P$ is generated by the $X_{i,j}$ (matrices with $B=C=0$) and the $Z_{i,j}$
(matrices with $A=B=D=0$). 
The images of these elements in $End(S^+)$ 
(as well as those of the $Y_{i,j}$) are:
{\renewcommand{\arraystretch}{1.5}
$$
so(2n)\,\longrightarrow\,End(S^+),\qquad
\left\{
\begin{array}{rclcl}
X_{i,j}&\longmapsto&
e_ie_{n+j}-\mbox{$\frac{1}{2}\delta_{ij}$}&\longmapsto&
L_{e_i}\circ D_{e_{n+j}}-\mbox{$\frac{1}{2}$}\delta_{ij},
\\
Y_{i,j}&\longmapsto&
e_ie_{j}&\longmapsto&
L_{e_i}\circ L_{e_j},
\\
Z_{i,j}&\longmapsto&
e_{i+n}e_{j+n}&\longmapsto&
D_{e_{i+n}}\circ D_{e_{j+n}}~.
\end{array}\right.
$$
} 

\

Since $D_{w^*}(1)=0$ for all $w^*\in W^*$, we see that $X_{i,j}$ and $Z_{i,j}$ map $1$
to an element in $\langle 1\rangle$. Hence $Lie(P)$ maps $\langle 1\rangle$ into itself and thus
also the inverse image of $P$ in $Spin(V)$ maps this line into itself.

The element $X_B\in so(2n)$ determined by $B$ is $X_B=\sum_{i<j} b_{ij}Y_{i,j}$.
It acts as left multiplication by $\omega_B:=\sum b_{ij}e_i\wedge e_j$ on $\wedge^{even}W$
and thus $exp(X_B)$ is left multiplication by $exp(\omega_B)\in \wedge^{even}W$. 
The exponential map of an endomorphism $\alpha$ is $\sum \alpha^n/n!$ and since the 
left multiplication 2-forms
generate a commutative subalgebra of nilpotent elements, $exp(X_B)$ is actually a finite sum
and one also has
$$
exp(\omega_B)\,=\,\prod_{i<j}exp(b_{ij}e_i\wedge e_j)\,=\,\prod_{i<j}(1+b_{ij}e_i\wedge e_j)~.
$$
We now show that, with $B_I$ the submatrix of $B$ with coefficients $b_{i,j}$ with $i,j\in I$,
$$
exp(\omega_B)\,=\,\sum_{I,\sharp I\,even}\,\mbox{Pfaff}(B_I)e_I~.
$$

In fact, $exp(\omega_B)\in \wedge^{even}W$ is a linear combination of the 
$e_I=e_{i_1}\wedge\ldots\wedge e_{i_{2k}}$, where  $i_1\leq\ldots\leq e_{i_{2k}}$
$I=\{i_1,\ldots,i_{2k}\}\subset\{1,\ldots,n\}$ is a subset with an even number of elements. 
Since for an integer $p$ one has that $\omega_B^p\in\wedge^{2p}W$,
the coefficient of $e_I$ is homogeneous of degree $2k$, with $2k=\sharp I$,
in the coefficients  $b_{ij}$ of $B$ and only those with $i,j\in I$
contribute. So the coefficient of $e_I$ is determined by the $2k\times 2k$ alternating 
submatrix $B_I$ of $B$ with rows and columns indexed by $I$. Moreover this coefficient is 
$(\sum_{i_k<i_l,i_k,i_l\in I} b_{i{_k}i{_l}}e_{i_k}\wedge e_{i_l})^k/k!$, which is indeed 
$\mbox{Pfaff}(B_I)$.

Since $\rho_V(\tilde{g}_B)Z_0\,=\,Z_B$ and $\gamma([Z_0])=1\in S^+$ 
we get $\gamma([Z_B])=\rho^+(\tilde{g})1=exp(\omega_B)\in S^+=\wedge^{even}W$.
The description of the spinor map follows immediately. For the equations defining the image
see \cite[III.3.2]{chevalley} or \cite{Lichtenstein}.
\qed

\

The following lemma is used several times in this paper, for example
to relate complex structures on $V_\RR$ to elements of $S^+_\CC$
or to weight two Hodge structures on $S^+$ as in the Kuga Satake construction. 
For $\dim V\neq 8$ however, $Spin(V)$ only allows one to
relate polarized weight two Hodge structures on $V$ to complex structures on the spin representation.
The special feature in the case $\dim V=8$ is triality, an automorphism of order three of $Spin(V)$,
which allows one to permute the three irreducible 8-dimensional representations $V,S^+,S^-$,
see \cite[\S 20.3]{FH}, \cite[Chapter 4]{chevalley}, and which is implicit in the proof of the lemma. 

\

\subsection{Lemma} \label{acth} 
\begin{enumerate}
 \item[a)]
Let $V=U\oplus U^*$ be a decomposition of $V=\CC^8$ with two maximally isotropic subspaces
with $[U],[U^*]\in IG(4,V_\CC)^+$.
For $t\in\CC$, $t\neq 0$, the orthogonal transformation $(t\mbox{id}_{U},t^{-1}\mbox{id}_{U^*})
\in (\mbox{End}(U)\oplus \mbox{End}(U^*))\cap SO(V)$ has a lift $h(t)\in Spin(V)$ 
which acts as follows on $S^+$:
$$
\rho^+(h(t))\,\ell_U\,=\,t^2 \ell_U,\quad \rho^+(h(t))\,\ell_{U^*}\,=\,t^2 \ell_{U^*},\qquad
\rho^+(h(t))\,s\,=\,s,\quad \forall s\in \langle \ell_U,\,\ell_{U^*}\rangle^\perp~,
$$
where $\ell_U,\ell_{U^*}\in S^+$ are (any) representatives of 
$\gamma([{U}]),\gamma([U^*])\in\PP S^+$.

\item[b)]
Let $Z_1,Z_2$ be two distinct maximally isotropic subspaces of $V_\CC$ in the family 
parametrized by $IG(4,V_\CC)^+$. Then $Z_1\cap Z_2=\{0\}$ if and only if
the complex plane $\langle [Z_1],[Z_2]\rangle$ is not contained in the spinor quadric $Q^+$.
\end{enumerate}

\ts 
We use that the spinor map is equivariant for the action of $Spin(V)$.
There is an element of $Spin(V)$ mapping $U$ to $W$ since $IG(4,8)^+=SO(V)/P$.
Then $U^*$ is mapped to $Z_B$ for some $B\in Alt_4$ and it is easy to see that there is another 
element in $Spin(V)$ fixing $W$ (so with $C=0$) and mapping $Z_B$ to $Z_0=W^*$.  
We thus may replace $W,W^*$ with $U,U^*$. 
The one parameter subgroup $h$ acts as multiplication by $t$ on $U\subset V$,
hence $h$ is generated by an $X\in {\mathfrak h}\subset 
spin(V)$ with $L_i(X)=1$ for $i=1,\ldots,4$ (and thus $X=\sum H_i$). 
The weights of $S^+$ are $(\pm L_1\pm L_2\pm L_3\pm L_4)/2$ with an
even number of $-$ signs, hence their values on $X$ are $2,-2$, with multiplicity one, 
and $0$ with multiplicity six. Thus $\rho^+(h(t))$ is semisimple with eigenvalues $t^2,t^{-2}$
and $1$, the last with multiplicity six. The eigenvalue $t^{-2}$, the lowest weight of $S^+$,
is on $Z_{U^*}$, see Theorem \ref{half-spinor}. 
The element $g\in SO(V)$ that maps $e_i\mapsto e_{i+4}$, $e_{i+4}\mapsto e_i$ for $i=1,\ldots,4$
interchanges $U$ and $U^*$ and acts (in the Adjoint representation) as $-id$ on $\mathfrak{h}$, 
hence the eigenvalue $t^2$ must be on $Z_U$. 
As $Spin(V)$ preserves $(\bullet,\bullet)_{S^+}$, the decomposition into these eigenspaces
is orthogonal. (For any $n$, the one parameter subgroup of $SO(V)$
that acts as multiplication by $t^2,t^{-2}$ on $e_1,e_{n+1}$ respectively and is the identity
on $\langle e_1,e_{n+1}\rangle^\perp$ is generated by an $X\in spin(V)$ with $L_1(X)=2$, $L_i(X)=0$
for $i\geq 2$ and thus $(1/2)(\pm L_1\pm L_2\ldots\pm L_n)(X)=\pm (1/2)X$, showing that 
the lift of this subgroup to $Spin(V)$ has only eigenvalues $t,t^{-1}$ on $S^+$,  
with the same multiplicities, and the same holds for $S^-$. 
A similar result holds for $SO(V)$ and its spin representation if $\dim V=2n+1$.)

Using the action of the orthogonal group,
if $Z_1\cap Z_2=\{0\}$, then we can map $Z_1,Z_2$ to $W,W^*$. As 
$[W]=e_*,[W^*]=1\in S^+$
and $(e_*,1)_{S^+}\neq 0$ it follows that the plane $\langle [Z_1],[Z_2]\rangle$
is not contained in $Q^+$. On the other hand, if $Z_1\cap Z_2\neq 0$, then we may assume $Z_1=W^*$ and
$Z_1=Z_B$ with $B$ the rank two alternating $4\times 4$ matrix with $\omega_B=e_1\wedge e_2$. Then
$[Z_1]=1$ and $[Z_2]=1+e_1\wedge e_2$ so that $\langle [Z_1],[Z_2]\rangle\subset Q^+$.
\qed

\

\subsection{The spinor map and the Pl\"ucker map}\label{soVreps}
We relate the spinor and Pl\"ucker maps on $IG(n,2n)^+$.
Even if the theory of line bundles on homogeneous spaces provides a natural setting for the results
below (cf.\ \cite[\S 23.3, p.393]{FH}, \cite[\S II]{BaloghHH}), we only use representation theory.

Let $\Gamma_\lambda$  be the 
irreducible $so(2n)$-representation with highest weight $\lambda$. 
The irreducible $so(2n)$-representation $S^\epsilon$, $\epsilon\in\{+,\,-\}$ 
has highest weight $\omega_n:=(L_1+\ldots+L_n)/2$, where $\epsilon=+$ if $n$ is even and $\epsilon=-$
else,
with highest weight vector $e_*:=e_1\wedge\ldots\wedge e_n\in S^\epsilon$, \cite[Prop.\  20.15]{FH}
(for $n=4$ we wrote $\alpha$ for $\omega_n$ in \S \ref{spinorpluecker}).

The highest weight of  $Sym^2(S^\epsilon)$ is thus $2\omega_n=L_1+\ldots+L_n$, 
with highest weight vector $e_*\odot e_*$.
In particular, $\Gamma_{2\omega_n}$ is an irreducible component of $Sym^2(S^+)$. 
In case $n=4$ we have $\dim \Gamma_{2\omega_n}=35=\dim Sym^2(\Gamma_{\omega_n})-1$ 
(see below for dimension formula) and thus 
(see \cite[Exercise 19.6]{FH} for general $n$):
$$
Sym^2(S^+)\,=\,Sym^2(\Gamma_{\omega_n})\,=\,\Gamma_{2\omega_n}\,\oplus\,\Gamma_0~,\qquad (n\,=\,4),
$$
where $\Gamma_0$ is the trivial 1-dimensional representation (for $n=4$ the representation
$\Gamma_{\omega_n}$ has 
an invariant quadratic form and thus is self-dual \cite[Exercise 20.38]{FH}), 
this quadratic form produces $\Gamma_0$).

Now we consider the representation of $so(2n)$ on $\wedge^nV$.
In the standard representation $V$ of $so(2n)$ the basis vector $e_i$ has weight $L_i$
(and the basis vector $e_{i+n}$ has weight $-L_i$). Thus 
$e_1\wedge\ldots\wedge e_n\in\wedge^nV$ has weight $L_1+\ldots+L_n=2\omega_n$ 
and it is the highest weight vector in $\wedge^nV$
(cf.\ \cite[\S 19.2]{FH}). Therefore $\Gamma_{2\omega_n}$ is a subrepresentation of $\wedge^nV$.
The $so(2n)$-representation $\wedge^nV$ is in fact reducible and it has two irreducible components 
of the same dimension (\cite[Remarks p.289-290; Exercise 24.43]{FH}), 
$$
\wedge^n V\,=\,\Gamma_{2\omega_n}\,\oplus\,\Gamma_{2\omega_{n-1}},\qquad 
\dim \Gamma_{2\omega_n}\,=\,\dim \Gamma_{2\omega_{n-1}}\,=\,\mbox{$\frac{1}{2}$}\binom{2n}{n}~,
$$
where $\omega_{n-1}:=(L_1+\ldots+L_{n-1}-L_n)/2$ is the highest weight of $S^{\epsilon'}$
(where $\{\epsilon,\epsilon'\}=\{+,-\}$).
This splitting is also obtained as the eigenspace decomposition for 
the Hodge star operator defined 
by the bilinear form $(\bullet,\bullet)_V$, cf.\ \cite[Example 3.5.2]{Joyce}. The two summands
are known as the selfdual and anti-selfdual forms.

The Pl\"ucker map, restricted to $IG(n,2n)^+$, is the natural map 
$\pi:IG(n,2n)^+\rightarrow \PP \wedge^nV$ and for $n$ even it is
thus the composition
$$
IG(n,2n)^+\,\stackrel{\gamma}{\longrightarrow}\,\PP S^+\,\stackrel{\nu}{\longrightarrow}\,
\PP\Gamma_{2\omega_n}\,\subset\,\PP\wedge^nV,
$$
where $\nu$ is the second Veronese map (for odd $n$ replace $\omega_n$ by
$\omega_{n-1}$, the highest weight of $S^+$).

Since the spinor map is given by Pfaffians and the Pl\"ucker map is given by minors 
on the open subset of $IG(n,2n)^+$ parametrized by alternating
matrices, this result implies
that any quadratic expression in Pfaffians is a linear combination of minors, see \cite{BaloghHH}.

\subsection{Cayley classes and $Spin(7)$}
We now restrict ourselves to the case $n=4$. 
The half-spin representation $\rho^+$ on $S^+$ maps the group $Spin(V)$ 
onto the orthogonal group $SO(S^+)$ 
of the bilinear form $(\bullet,\bullet)_{S^+}$.
For any $s\in S^+$ with $(s,s)_{S^+}\neq 0$, the stabilizer of $s$ in $SO(S^+)$ 
is the orthogonal group $SO(s^\perp)\cong SO(7)$.

The inverse image of $SO(s^\perp)$
in $Spin(V)$ is denoted by $Spin(V)_s$ and it is isomorphic to $Spin(7)$.
In the standard representation $\rho_V$ of $Spin(V)$ on $V$,
the subgroup $Spin(V)_s$ still acts irreducibly, in fact $V$ is isomorphic with the (unique,
irreducible) spin representation of $Spin(7)$.

\subsection{Representations of $spin(V)_s=so(7)$}\label{spin7}
Recall from \S \ref{cayleyclasses} that the image of
$s\odot s$ under the composition 
$Sym^2(S^+)\rightarrow \Gamma_{2\alpha}\hookrightarrow \wedge^4V$ is called the Cayley class
$c_s$ of $s$. Since $s$ is fixed by $Spin(V)_s$, the 4-form $c_s$ is also fixed by the 
Lie algebra $spin(V)_s\cong so(7)$. 
We now show that $c_s$ is the unique $spin(V)_s$-invariant in $\wedge^4V$
by considering the restriction to $so(7)$ of the $so(V)=so(8)$-representations considered
in \S \ref{soVreps}.

Multiplication by $s$ gives an inclusion of $spin(V)_s$-representations
$$
\begin{array}{rcl}
S^+\,=\,\langle s\rangle\,\oplus\,s^\perp\,\hookrightarrow\,Sym^2(S^+)&=& 
\Gamma_0\,\oplus\,\Gamma_{2\omega_n}\\
&=&\Gamma_0\,\oplus\,\langle c_s\rangle
\,\oplus\,s\odot s^\perp\,\oplus\,\Gamma_{(2,0,0)} \\
&=&\Gamma_0\,\oplus\,\Gamma_{(0,0,0)}\,\oplus\,\Gamma_{(1,0,0)}\,\oplus\,\Gamma_{(2,0,0)}~,
\end{array}
$$
where $\Gamma_0$ and $\Gamma_{(0,0,0)}$ are trivial $spin(V)_s$-representations,
$\Gamma_{(1,0,0)}\cong s\odot s^\perp\cong s^\perp$ is the standard seven dimensional 
representation of $spin(V)_s\cong so(7)$ and $\Gamma_{(2,0,0)}$ is irreducible of 
dimension $35-1-7=27$ (the notation $\Gamma_{(a,b,c)}$ for $so(7)$-representations
is as in \cite{FH}). 

The representation of $spin(V)_s$ on the $spin(V)$-representation $\Gamma_{2\omega_n}$ 
is thus a direct sum of three irreducible representations. Its representation 
on the other irreducible component $\Gamma_{2\omega_{n-1}}$ of $\wedge^4V$
is irreducible and it is isomorphic to $\Gamma_{(0,0,2)}$. Thus one has the 
$spin(7)=so(7)$-decomposition into irreducible representations (cf.\ \cite[Prop 2]{Munoz},
\cite[Prop.\ 10.5.4]{Joyce}):
$$
\wedge^4V\,=\,\Gamma_{(0,0,0)}\,\oplus\,\Gamma_{(1,0,0)}\,\oplus\,\Gamma_{(2,0,0)}\,
\oplus\,\Gamma_{(0,0,2)}~,
$$
since there is a unique copy of the trivial representation of $so(7)$ in $\wedge^4V$,
the Cayley class is the unique $spin(V)_s$ invariant in $\wedge^4V$.

\

\subsection{}
The following proposition computes the 4-form $c_s$, which spans the trivial 
$spin(V)_s$-subrepresentation
$\Gamma_{(0,0,0)}$ in $\wedge^4V$,
explicitly in a case of interest in 
Markman's paper, cf.\ \cite[1.4.1, Proposition 11.2]{Markman}.
There $s$ is called $w=s_n$. 
We consider in fact $\frac{1}{n+1}c_w$ and 
we write $n$ for his $n+1$. 
Notice that the computation below uses only representation theory.

\subsection{Proposition}\label{explicitCs}
Let $n\in\ZZ$, $n\neq 0$, and let $s=s_n=1-n e_* \in S^+$.
where $e_*:=e_1\wedge e_2\wedge e_3\wedge e_4\in\wedge^{even} W=S^+$.
Then we have, up to a scalar multiple,
$$
c_s\,=\,-n\alpha^2\,+\,4n^2\beta\,+4\gamma\quad(\in \wedge^4V)~,
$$
where the forms, now in $\wedge^*V$, involved are:
$$
\alpha\,:=\,e_1\wedge e_5+\ldots+e_4\wedge e_8,\qquad
\beta\,:=\,e_1\wedge \ldots\wedge e_4,\qquad \gamma\,:=\,e_5\wedge\ldots\wedge e_8~.
$$

\ts 
The space of $spin(V)_s$-invariants in $\wedge^4V$ is one dimensional and it is spanned by $c_s$,
see \S \ref{spin7}. 
So it suffices to show that the right hand side is a non-zero $spin(V)_s$-invariant form.

The Lie algebra $spin(V)_{1,e_*}$
that acts trivially on the two dimensional subspace of $S^+$ spanned by $1,e_*$
is isomorphic to $so(6)\cong sl(4)$. The representation of $sl(4)$ on $V=W\oplus W^*$ is reducible
and $W$ is the standard representation of $sl(4)$ whereas $W^*$ is the dual of the standard 
representation. This implies that $\beta\in \wedge^4W\subset \wedge^4V$ and $\gamma\in \wedge^4W^*
\subset \wedge^4V$ as well as the 2-form $\alpha$, which is the $sl(4)$-invariant in 
$W\otimes W^*\subset \wedge^2V$ corresponding to the symplectic form
$((w_1,w_1^*),(w_2,w_2^*))=w_1^*(w_2)-w_2^*(w_1)$ on $V$, are $spin(V)_{1,e_*}$-invariants.
On the other hand, 
$$
\wedge^4(W\oplus W^*)\,=\,\wedge^4W\,\oplus W\otimes \wedge^3W^*\,\oplus\,
\wedge^2W\otimes \wedge^2W^*\,\oplus\,\wedge^3W\otimes W^*\,\oplus\,\wedge^4W^*~.
$$
Since $W,W^*$ have dimension four, $\wedge^3W^*\cong W$ and it is well-known that there
are no $sl(4)$-invariants in $W\otimes W$ nor in $W^*\otimes W^*$. Also $\wedge^2W$ is irreducible
and thus the $sl(4)$-invariants in $\wedge^2W\otimes \wedge^2W^*\cong End(\wedge^2W)$ 
are a one dimensional subspace spanned by the trace. Hence the subspace of $sl(4)$-invariants in
$\wedge^4V$ has dimension three. Since $\alpha^2,\beta,\gamma$ are linearly independent
invariants, the invariant subspace is
$$
(\wedge^4 V)^{spin(V)_{1,e_*}}\,=\,(\wedge^4 V)^{sl(4)}\,=\,\langle\alpha^2,\,\beta,\gamma\rangle~.
$$

Since $spin(V)_{1,e_*}\subset spin(V)_s$ it remains to show that $c_s$ is the linear
combination of $\alpha^2,\beta,\gamma$ that is $spin(V)_s$ invariant.
This 21-dimensional Lie algebra is defined by
$$
spin(V)_s\,=\,\{X\,\in\,spin(V):\;Xs\,=\,0\,\}
$$
and the action of $spin(V)$ on $S^+$ is given in the proof
of Theorem \ref{half-spinor}. It is then easy to check that the following elements 
(of $so(2n)\cong spin(V)$)
span $spin(V)_s$:
$$
{\mathfrak h}_s\,:=\,\{\sum a_iX_{i,i}:\;\sum a_i\,=\,0\,\},\quad 
X_{i,j}\quad (i\neq j),\quad nY_{i,j}\,\pm\,Z_{k,l}\quad (\{i,j,k,l\}=\{1,\ldots,4\}),
$$
where the sign depends on $i,\ldots,l$. In particular, $X:=nY_{1,2}+Z_{3,4}\in spin(V)_s$
(in fact $X$ acts as $e_1e_2+D_{e_3}D_{e_4}$ on $S^+$ and $X(1)=ne_1e_2$, 
$X(e_*)=-e_1e_2$, so $Xs=0$).
The action of $X$ on $V$ is given by
$$
\begin{array}{rclrclrclrcl}
                                X(e_1)&=&0,\;&X(e_2)&=&0,\;&X(e_3)&=&-e_8,\;&X(e_4)&=&e_7,\\
                               X(e_5)&=&-ne_7,\;&X(e_6)&=&ne_8,\;&X(e_7)&=&0,\;&X(e_8)&=&0~.
                               \end{array}
$$
Since the Lie algebra element $X$ acts a derivation on $\wedge^4V$ we have
$$
X(\alpha)\,=\,X(e_1)\wedge e_5+e_1\wedge X(e_5)+\ldots\,=\,-2ne_1\wedge e_2\,+\,2e_7\wedge e_8~.
$$
Thus 
$$
X(\alpha^2)\,=\,2\alpha\wedge X(\alpha)\,=\,-4n(e_1\wedge e_2)\wedge(e_3\wedge e_7+e_4\wedge e_8)
+4(e_1\wedge e_5+e_2\wedge e_6)\wedge (e_7\wedge e_8)~.
$$
Similarly one finds
$$
X(\beta)\,=\,(e_1\wedge e_2)\wedge (e_4\wedge e_8\,+\,e_3\wedge e_7),\qquad 
X(\gamma)\,=\,-n(e_2\wedge e_6\,+\,e_1\wedge e_5)\wedge (e_7\wedge e_8)~.
$$
Therefore the only non-trivial linear combination of $\alpha^2,\beta,\gamma$ that is mapped to zero 
by $X$ is $-n\alpha^2+4n^2\beta+4\gamma$. Hence this must be the unique $spin(V)_s$-invariant in
$\wedge^4V$.
\qed

\

\subsection{Kuga Satake varieties}\label{KSvar}
Let $S^+$ be the lattice introduced in \S \ref{spinormap} (and not its complexification).
As in Theorem \ref{Weil}, let $h,s\in S^+\cong U^{\oplus 4}$ be two perpendicular elements such that their span is 
a positive definite sublattice. 
Let $H=H_{h,s}$ be the rank $6$ sublattice
of signature $(2+,4-)$ orthogonal to $\langle h,s\rangle$:
$$
H\,:=\,\langle h,s\rangle^\perp\,=\,\{t\in S^+:\,(t,h)\,=\,(t,s)\,=\,0\,\}~.
$$
With this notation we have
$$
\Omega_{\{ h,s\}^\perp}\,=\,\{\ell\in\PP H_\CC:\quad 
(\ell,\ell)_{S^+}\,=\,0,\quad (\ell,\bar{\ell})_{S^+}\,>\,0\;\}~.
$$

Recall that any $\ell\in \Omega_{\{ h,s\}^\perp}$
defines an abelian fourfold of Weil type with underlying torus $\cT_\ell$
by Theorem \ref{Weil}.
Such an $\ell$ also defines a weight two Hodge structure on $H$ denoted by $H_\ell$
as follows:
$$
H_{\ell,\CC}\,=\,H_\CC\,=\,\oplus_{p+q=2} H^{p,q}_\ell,\qquad
H^{2,0}_\ell\,:=\,\CC\ell,\qquad H^{0,2}_\ell\,:=\,\CC\bar{\ell},\qquad H^{1,1}_\ell\,=\,
\big(H^{2,0}_\ell\oplus H^{0,2}_\ell\big)^\perp~.
$$
This Hodge structure is polarized since the restriction of $(\bullet,\bullet)_{S^+}$ to 
the two dimensional real subspace $(H^{2,0}_\ell\oplus H^{0,2}_\ell)\cap H_\RR$ is positive definite.

As $\dim H^{2,0}_\ell=1$, there is a Kuga Satake (abelian) variety $A_\ell$, of dimension $16$,
associated to $H_\ell$ (see \cite{KugaS}, \cite{Deligne}, \cite{vG_KS}).
In general, it has the property that $H_\ell$ is a Hodge substructure of $H^2(A_\ell^2,\QQ)$,
but in this case there are actually several copies of $H_\ell$ in $H^2(A_\ell,\QQ)$,
see \S \ref{H2ell}.
The even Clifford algebra $C(H)^+$ of $H$ is a lattice in the 
real vector space $C(H)^+\otimes_\ZZ\RR$ of dimension $2^5=32$. 
A complex structure on $C(H)^+_\RR$ is defined by left multiplication
by $f_1f_2\in C(H)^+_\RR$, with $f_1,f_2\in H_\RR$ such that
$(f_1,f_1)_{S^+}=1$ and $H^{2,0}_\ell=\langle f_1+if_2\rangle$ 
(cf.\ \cite[\S 5.6]{vG_KS}). 
The abelian variety $A_\ell$ is the quotient $(C(H)^+_\RR,f_1f_2)/C(H)^+$.

In \cite[Cor.\ 6.3, Thm 6.4]{Lombardo} it is shown that $A_\ell$ is isogeneous to 
$B_\ell^4$, where $B_\ell$ is an abelian fourfold of Weil type with 
trivial discriminant. 
The following proposition, due to O'Grady (\cite[\S 5.3]{O'G}), 
shows that $B_\ell$ and $\cT_\ell$ are isogeneous.
In \cite{O'G} one finds a more explicit description of this result, 
as well as applications to generalized Kummer varieties.

\subsection{Proposition} \label{cTKS}
For $\ell\in \Omega_{\{ h,s\}^\perp}$ 
the Kuga Satake variety $A_\ell$ of the polarized weight 
two Hodge structure $H_\ell$ is isogeneous to $\cT_\ell^4$, 
where $\cT_\ell$ is the abelian fourfold of Weil type defined by $\ell$. 

\ts
The right multiplication on $C(H)^+_\RR$ by an element of $C(H)^+$ 
preserves the lattice, commutes with the complex structure and thus defines an element in 
$\mbox{End}(A_\ell)$.
The $\QQ$ vector space $H_\QQ$ is not a direct sum of two maximally isotropic subspaces
and,
whereas $C(H)^+_\CC\cong M_4(\CC)\oplus M_4(\CC) $ (as in \S \ref{cliffordV}),
one now has an isomorphism of algebras (\cite[Thm.\ 6.2]{Lombardo}), 
where $M_4(K)$ are the $4\times 4$ matrices with coefficients in $K$,
$$
C(H)^+_\QQ\,:=\,C(H)^+\otimes_\ZZ\QQ\,\cong\,M_4(K)\,\subseteq \,
\mbox{End}(A_\ell)_\QQ,\qquad K\,:=\,\QQ(\sqrt{-ab}).
$$
This implies that any $A_\ell$ is isogeneous to $B_\ell^4$, 
where $B_\ell$ is an abelian fourfold with $K\subset \mbox{End}(B_\ell)_\QQ$  
($B_\ell$ is only determined up to isogeny). 

It remains to show that $B_\ell$ and $\cT_\ell$ are isogeneous.
The inclusion $Spin(H)\subset Spin(S^+)=Spin(V)$ defines a representation of $Spin(H)$ on $V$
which is its spin representation. The isomorphism $C(H)^+_\QQ\cong M_4(K)$ implies that 
$$
C(H)^+_\QQ\,\cong \,V_\QQ^{\oplus 4}
$$ 
as $Spin(H)$-representations. The same holds with $\QQ$ replaced by $\RR$.
The weight two Hodge structure on the $Spin(H)$-representation $H_\ell$ is defined
by the one parameter subgroup $h_\ell$ of $Spin(H)_\RR\subset Spin(S^+)_\RR$
introduced in the proof of Proposition \ref{CayleyHodge}. 
In fact,
$h_\ell(t)$ acts on $S^+$ as multiplication by $t^2$ on $\CC\ell$, by $t^{-2}$ on $\CC\bar{\ell}$
and it is trivial on $\langle\ell,\bar{\ell}\rangle^\perp$.
The complex structure on $C(H)^+_\RR\cong V_\RR^{\oplus 4}$, which defines 
the Kuga Satake variety $A_\ell\sim B_\ell^{\oplus 4}$, 
is also defined by $h_\ell$ (\cite[Prop.\ 6.3]{vG_KS}), now acting on $V_\RR^4$. 
As $\rho_V(h_\ell)=h_{V,\ell}$, the complex structure is $J_\ell$ on $V_\RR$.
It follows that $B_\ell$ and $\cT_\ell$ are isogeneous.
\qed

\

\subsection{Remark} \label{remark_spin6}
The proof of Proposition \ref{cTKS} uses the (algebraic) subgroup $Spin(H)=Spin_{h,s}$ of
$Spin(S^+)=Spin(V)$. The decomposition $S^+_\QQ=H_\QQ\oplus R_\QQ$, with $R:=\langle h,s\rangle$,
implies that we actually have two commuting subgroups $Spin(H), Spin(R)\subset Spin(S^+)$.

Recall from \S \ref{endomorphisms} that $R_\CC=\CC\kappa\oplus\CC\bar{\kappa}$ with 
$\kappa,\bar{\kappa}\in Q^+$.
The decomposition of $V_\CC=Z_{\kappa,\CC}\oplus
Z_{\bar{\kappa},\CC}$ in the two isotropic eigenspaces for the $K$-action defines,
as in Lemma \ref{acth}, a one parameter subgroup $h_R$ of $Spin(S^+_\RR)$.
As $h_R(t)\kappa=t^2\kappa$, $h_R(t)\overline{\kappa}=t^{-2}\bar{\kappa}$, this identifies
the subgroup $Spin(R_\RR)$ with this one parameter subgroup, $h_R(U(1))=Spin(R_\RR)$.
In particular, the $K$-action on $V_\QQ$ is generated by $Spin(R)$ 
and the scalar multiples of the identity.

The fact that $Spin(H), Spin(R)\subset Spin(S^+)$ commute implies that the subspaces 
$Z_{\kappa,\CC},Z_{\bar{\kappa},\CC}$ are $Spin(H_\CC)$-invariant subspaces. Thus the 
spin representation of $Spin(H_\CC)$ on $V_\CC$ is reducible. 
These two subspaces are the two half-spin representations of $Spin(H_\CC)$. 

There is an isomorphism $Spin(H_\CC)\cong SL(4,\CC)$
and the half-spin representations are identified with the standard representation $\CC^4$ of
$SL(4,\CC)$ and its dual $(\CC^4)^*$. The representation $H_\CC$ is identified with 
$\wedge^2 \CC^4\cong \wedge^2 (\CC^4)^*$, the isomorphism follows from the pairing,
defined by the wedge product,
$(\wedge^2\CC^4)\times(\wedge^2\CC^4)\rightarrow\wedge^4\CC^4\cong\CC$.

\

\subsection{The second cohomology group of $\cT_\ell$} \label{H2ell}
In \cite{Lombardo} the 
Hodge structure on the second cohomology group  $H^2(B,\QQ)$
of an abelian fourfold of Weil type with field $K$ is studied. 
This group has dimension $\binom{8}{2}=28$
and decomposes under the $K$-action into a $16=1+15$-dimensional subspace $S_B'$ 
on which  $x\in K$ acts 
as $x\bar{x}$, this subspace includes the polarization of Weil type. 
There is a complementary subspace $S_B$
on which the eigenvalues of $x$ are $x^2,\bar{x}^2$ of dimension $12$. This subspace 
can be identified with the six dimensional $K$ vector space $\wedge_K^2H^1(B,K)$. 
$$
H^2(B,\QQ)\,=\,S_B\,\oplus\,S'_B,\qquad 
S_B'\,:=\,\{\xi\in H^2(B,\QQ):\,x^* \xi\,=\,x\bar{x}\xi,\quad\forall x\in K\,\}~.
$$
For a general fourfold of Weil type 
(so $SMT(B)_\RR\cong SU(2,2)$) the Hodge structure $S_B$ is a simple Hodge structure (so does not
admit non-trivial Hodge substructures) if and only if the discriminant of $B$ is non-trivial
\cite[Cor.\ 3.6]{Lombardo}. 

In case the discriminant is trivial, one finds that 
$S_B\cong H_B^{\oplus 2}$, for a weight two, rank six, polarized, Hodge structure $H_B$ 
which has  Hodge numbers $(1,4,1)$. Moreover, the Kuga Satake variety of $H_B$ is isogeneous to 
$B^4$, so one recovers the weight two Hodge structure $H_B$ from its Kuga Satake variety.

The following proposition uses this result to show that the abelian fourfolds of Weil type $\cT_\ell$
have trivial discriminant.

\subsection{Proposition}\label{proptrivdisc}
For $\ell\in \Omega_{\{h,s\}^\perp}$, with $h,s$  as in Theorem \ref{Weil},
the polarized abelian fourfold of Weil type $(\cT_\ell,K,\omega_K)$ has trivial discriminant.

\ts
By \cite[Cor.\ 3.6]{Lombardo} it suffices to show that 
$(H^2_{\ell,\QQ})^{\oplus 2}$ is isomorphic to the Hodge substructure 
$S_{\cT_\ell}\subset H^2(\cT_\ell,\QQ)$.

As in the proof of Proposition \ref{CayleyHodge}, 
the (weight one) Hodge structure on $V=H^1(\cT_\ell,\ZZ)$ 
defines a one parameter subgroup $h_\ell$ in $Spin(V)$ 
(actually in $Spin(V)_{h,s}\subset Spin(S^+)=Spin(V)$).
A representation $U$ of $Spin(V_\RR)$ on a real vector space $U$ defines a Hodge decomposition
$U_\CC=\oplus U^{p,q}$, with $\overline{U^{p,q}}=U^{q,p}$,
given by the eigenspaces $U^{p,q}=\{u\in U:h_\ell(z)u=z^a\bar{z}^bu$ 
(but the weight is not uniquely defined since $z\bar{z}=1$).

The representation $\rho^+$ on $S^+_\RR$ has the Hodge decomposition
$$
(S^+)^{2,0}\,=\,H^{2,0}_\ell\,=\,\CC\ell,
\qquad (S^+)^{0,2}\,=\,\overline{(S^+)^{2,0}},
\qquad \qquad (S^+)^{1,1}\,=\,\big((S^+)^{2,0}\oplus (S^+)^{0,2}\big)^\perp
$$
since these spaces are the eigenspaces for $h_\ell$ acting on $S^+_\CC$ (see Lemma \ref{acth}).
The Hodge structure $S^+_\QQ$ is a direct sum of Hodge structures
$$
S^+_\QQ\,=\,H^2_{\ell,\QQ}\,\oplus\,R_\QQ,\qquad R\,:=\,\langle h,s\rangle~,
$$
where
$R_\QQ\cong\QQ(-1)^2$ is a trivial Hodge substructure with $R_\QQ^{1,1}=R_\CC$.

There is an isomorphism of $Spin(V)=Spin(S^+)$-representations $\wedge^2S^+=\wedge^2 V$
(both are the irreducible $so(8)$-representation with highest weight 
$(L_1+L_2+L_3+L_4)/2+(L_1+L_2-L_3-L_4)/2=L_1+L_2$).
Hence we get a splitting
of the Hodge structure on $\wedge^2 S^+_\QQ$ (which is again defined by $h_\ell$ eigenspaces)
in three Hodge substructures which have dimensions 
$\binom{6}{2}=15$, $6\cdot 2=12$ and $1$ respectively:
$$
\wedge^2 S^+_\QQ\,=\,
(\wedge^2 H^2_{\ell,\QQ})\;\oplus\; (H^2_{\ell,\QQ}\otimes R_\QQ)\;\oplus\;(\wedge^2R_\QQ)~.
$$
(The Hodge structure $S^+$ has weight two, so the Hodge structure on $\wedge^2S^+$ 
should have weight four.
However, $(\dim S^+)^{2,0}=1$,
so $\wedge^2S^+_\QQ$ has trivial $(4,0)$ and $(0,4)$ summands 
and thus it is the Tate twist of a weight two Hodge structure.) 

Using the isomorphisms $\wedge^2 S^+_\QQ=\wedge^2V=H^2(\cT_\ell,\ZZ)$ we see that
$$
H^2_{\ell,\QQ}\otimes R_\QQ\,\cong\, (H^2_{\ell,\QQ})^{\oplus 2}\;\hookrightarrow\;H^2(\cT_\ell,\QQ)
$$ 
is a non-simple Hodge substructure of $H^2(\cT_\ell,\QQ)$. 

It remains to check that $x\in K$ has eigenvalues $x^2,\bar{x}^2$ on this substructure. 
One can deduce this from the fact that 
representation $\wedge^2V_\CC$ of the complex Mumford Tate group $SL(4,\CC)$  of $\cT_\ell$
is isomorphic to 
$$
\wedge^2(\CC^4\oplus(\CC^4)^*)\cong (\wedge^2\CC^4)^{\oplus 2}\,\oplus\, \CC^4\otimes(\CC^4)^*
$$
and the last summand is the direct sum of a trivial one dimensional representation and an irreducible
15 dimensional representation. As the complexification of a Hodge substructure is 
a subrepresentation, there is a unique subrepresentation of dimension $12$.
Hence $S_{\cT_\ell}=H^2_{\ell,\QQ}\otimes R_\QQ$ as desired.

Alternatively, by Remark \ref{remark_spin6}, the $K^\times$-action
is essentially given by the subgroup $Spin(R)$ of $Spin(S^+)$. This subgroup acts trivially
on $\wedge^2H^2_{\ell,\QQ}$ and $\wedge^2R_\QQ$, so $K$
acts through the norm on these summands. Therefore $S_{\cT_\ell}=H^2_{\ell,\QQ}\otimes R_\QQ$
is not simple.
\qed

\quad

\

\end{document}